\documentclass[12pt]{article}

\usepackage{amssymb}
\usepackage{amsmath}
\newtheorem{claim}{\bf \t}[part]

\numberwithin{Assumption}{section} \numberwithin{Corollary}{section}
\numberwithin{Definition}{section} \numberwithin{equation}{section}
\numberwithin{Example}{section} \numberwithin{Lemma}{section}
\numberwithin{Proposition}{section} \numberwithin{Remark}{section}
\numberwithin{Theorem}{section}

\def\v{\varepsilon}

\def\x{\xi}
\def\t{\theta}

\def\mb{\mathbf}
\def\a{\alpha}

\def\g{\gamma}
\def\d{\delta}
\def\l{\lambda}
\def\f{\frac}

\def\r{\rho}
\def\ra{\rightarrow}

\def\di{\displaystyle}
\def\i{\infty}

\def\text#1{{\rm #1}}

\begin{document}
\date{}
\title{\Large \bf Stability of Rarefaction Waves to the 1D Compressible Navier-Stokes Equations with Density-dependent Viscosity}

\author{\small \textbf{Quansen Jiu},$^{1,3}$\thanks{The research is partially
supported by National Natural Sciences Foundation of China (No.
10871133) and Project of Beijing Education Committee. E-mail:
qsjiumath@gmail.com}
 \qquad  \textbf{Yi Wang}$^{2,3}$\thanks{The research is partially
supported by National Natural Sciences Foundation of China (No.
10801128). E-mail: wangyi@amss.ac.cn. }\qquad and \qquad
\textbf{Zhouping Xin}$^{3}$\thanks{The research is partially
supported by Zheng Ge Ru Funds, Hong Kong RGC Earmarked Research
Grant CUHK4042/08P and CUHK4040/06P, and a Focus Area Grant at The
Chinese University of Hong Kong. Email: zpxin@ims.cuhk.edu.hk}
\\
\small $^1$School of Mathematical Sciences, Capital Normal
University, Beijing 100048, P.R. China \\\small $^2$Institute of
Applied Mathematics, AMSS, CAS, Beijing 100190, P. R. China
 \\\small $^3$The Institute of Mathematical Sciences, CUHK, Shatin N.T., HongKong}

\small \maketitle { \bf Abstract:} In this paper, we study the
asymptotic stability of rarefaction waves for the compressible
isentropic Navier-Stokes equations with density-dependent
viscosity. First, a weak solution around a rarefaction wave to the
Cauchy problem is constructed by approximating the system and
regularizing the initial values which may contain vacuum state.
Then some global in time estimates on the weak solution are
obtained. Based on these uniform estimates, the vacuum states are
shown to vanish in finite time and the weak solution we
constructed becomes a unique strong one. Consequently, the
stability of the rarefaction wave is proved in a weak sense. The
theory holds for large-amplitudes rarefaction waves and arbitrary
initial perturbations.


 {\bf AMS SC2000:} 35L60, 35L65

\section{Introduction } \setcounter{equation}{0}
\setcounter{Assumption}{0} \setcounter{Theorem}{0}
\setcounter{Proposition}{0} \setcounter{Corollary}{0}
\setcounter{Lemma}{0} In this paper, we consider the following
compressible and isentropic Navier-Stokes equations with
density-dependent viscosity
\begin{eqnarray}\label{(1.1)}
\begin{cases}
    \rho   _{t}+ (\rho u) _{x} = 0,   \qquad\qquad\qquad\qquad\qquad\qquad x\in\mathbf{R},t>0,\cr
    (\rho u)_{t}+  \big {(} \rho u^2 +p(\r)\big {)}_x
    =(\mu(\r)  u_{x})_x,\cr
   \end{cases}
\end{eqnarray}
where $\rho (t, x)\geq 0 $, $u(t, x) $ represent the density and the
velocity of the gas, respectively. Assume that the pressure and
viscosity function is given by
\begin{eqnarray}
p(\r)=A \r^{\gamma },\qquad \mu(\r)=B\r^\a,~~\label{(1.2)}
\end{eqnarray}
where $\gamma >1$ denotes the adiabatic exponent, $\a>0$  and
$A,B>0$ are the gas constants. Without loss of generality, we assume
that $A=B=1$.

We consider the Cauchy problem (\ref{(1.1)}) with the initial values
\begin{equation}
  (\r,\r u)( 0, x)=(\r_0,m_0)(x)\rightarrow (\r_{\pm},m_{\pm}) ~~ \text{as} ~~ x \rightarrow
 \pm\infty,\label{(1.3)}
\end{equation}
where $\r_{\pm}$,  $m_{\pm} $  are given constants. Here we assume
that there is no vacuum at the far field, i.e., $\r_{\pm}>0$, thus
we can define the far field velocity by $u_\pm=\f{m_\pm}{\r_\pm}$.

The large time asymptotic behavior of solutions to (\ref{(1.1)})
is expected to be closely related to that of the corresponding
Euler system
\begin{equation}
\left\{
\begin{array}{l}
\di \r_t+(\r u)_x=0,\\
\di (\r u)_t+(\r u^2+p(\r))_x=0,\\
\end{array}
\right. \label{(1.4)}
\end{equation}
 The Euler system (\ref{(1.4)}) is a strictly hyperbolic one for
 $\r>0$ whose characteristic fields are both genuinely nonlinear, that
 is, in the equivalent system
$$
\left(
\begin{array}{l}
\di \r\\
\di u
\end{array}
\right)_t + \left(
\begin{array}{cc}
\di u&\quad \r\\
\di p^\prime(\r)/\r&\quad u
\end{array}\right)\left(
\begin{array}{l}
\di \r\\
\di u
\end{array}\right)_x=0,
$$
the Jacobi matrix
$$
\left(
\begin{array}{cc}
\di u&\quad \r\\
\di p^\prime(\r)/\r&\quad u
\end{array}\right)
$$
has two distinct eigenvalues
$$
\l_1(\r,u)=u-\sqrt{p^\prime(\r)},\qquad
\l_2(\r,u)=u+\sqrt{p^\prime(\r)}
$$
with corresponding right eigenvectors
$$
r_i(\r,u)=(1,(-1)^i\f{\sqrt{p^\prime(\r)}}{\r})^t,\qquad i=1,2,
$$
such that
$$
r_i(\r,u)\cdot \nabla_{\r,u}\l_i(\r,u)=(-1)^i\f{\r
p^{\prime\prime}(\r)+2p^\prime(\r)}{2\r\sqrt{p^\prime(\r)}}\neq
0,\quad i=1,2.
$$
Define the $i-$Riemann invariant $(i=1,2)$ by
$$
\Sigma_i(\r,u)=u+(-1)^i\int^\r\f{\sqrt{p^\prime(s)}}{s}ds
$$
such that
$$
\nabla_{(\r,u)}\Sigma_i(\r,u)\cdot r_i(\r,u)\equiv0,\qquad \forall
\r>0,u.
$$
There are two families of rarefaction waves to the Euler system
(\ref{(1.4)}). Here we only consider $2-$rarefaction wave, which is
characterized by the fact that $2-$Riemann invariant
$\Sigma_2(\r,u)$ is constant in $(x,t)$ and $2-$characteristic
speed, $\l_2(\r,u)$ is increasing in $x$. Suppose that the end
states $(\r_\pm,u_\pm)$ of the initial values of \eqref{(1.4)}
satisfy
$$
\Sigma_2(\r_-,u_-)=\Sigma_2(\r_+,u_+),\qquad\l_2(\r_+,u_+)>\l_2(\r_-,u_-),
$$
Then the state $(\r_-,u_-)$ can be connected to the state
$(\r_+,u_+)$ by a $2-$rarefaction wave.  The $2-$rarefaction waves
for \eqref{(1.4)} connecting $(\r_-,u_-)$ to $(\r_+,u_+)$ converge
to each other time-asymptotically in supreme norm. For definiteness,
we choose a particular 2-rarefaction wave $(\bar\r,\bar u)(t,x)$ of
\eqref{(1.4)} defined by \eqref{(2.3)} in the next section.

In the present paper, we are interested in the large
time-asymptotic stability of the above rarefaction wave to the
density-dependent Navier-Stokes equation \eqref{(1.1)} with large
initial values which may contain vacuum states.

There have been large literature on the global existence and the
large time behavior of solutions to the system (1.1) and even in
the higher dimensional case when the viscosity $\mu(\r)$ is
constant. We refer to  \cite{FNP00}, \cite{JXZ}, \cite{KS},
\cite{L98}, \cite{Liu-Xin-1}, \cite{Matsumura-Nishihara-1},
\cite{Matsumura-Nishihara-2} and the references therein. However,
the possible appearance of the vacuum is one of the major
difficulties to the global existence of the solutions to the
compressible Navier-Stokes equations with constant viscosity. Hoff
and Smoller \cite{Hoff-Smoller} proved that weak solutions of the
compressible Navier-Stokes equations in 1D case do not contain
vacuum states in finite time if there is no vacuum initially. On
the other hand, Xin \cite{Xin} proved that there is no global
smooth solution to the Cauchy problem to compressible
Navier-Stokes equations with a nontrivial compactly supported
initial density, which means that the solution will blow up in
finite time in the presence of the vacuum states. Moreover, Hoff
\cite{Hoff} proved the results of the failure of the continuous
dependence of the weak solutions containing vacuum states on the
initial values.

Thus, when the solution may contain vacuum states, it is natural
to consider the compressible Navier-Stokes equations with
density-dependent viscosity, as was derived from Chapman-Enskog
expansions from the Boltzmann equation where the viscosity depends
on the temperature and thus on the density for isentropic flows.
Moreover, a one-dimensional viscous Saint-Venant system for the
shallow water, derived rigorously from the incompressible
Navier-Stokes equation with a moving free surface by
Gerbeau-Perthame recently in \cite{GP}, is expressed exactly as in
\eqref{(1.1)}-\eqref{(1.2)} with $\a=1$ and $\g=2$.

There are  many literatures on mathematical studies on
\eqref{(1.1)}-\eqref{(1.2)}. If the initial density is assumed to
be connected to vacuum with discontinuities, Liu, Xin and Yang
first obtained in \cite{LXY} the local existence of weak
solutions.  The global existence of weak solutions was obtained
later by \cite{J}, \cite{JXZ}, \cite{OMM}, \cite{YYZ}
respectively. If the initial density connects to vacuum
continuously, then new difficulty is encountered since no positive
lower bound for the density is available. This case is studied by
\cite{FZ},\cite{VYZ},\cite{YYZ} and \cite{YZ2} respectively.
However, most of these results concern with free boundary
problems. Recently, initial-boundary-value problems for
one-dimensional equations \eqref{(1.1)}-\eqref{(1.2)} with
$\mu(\rho)=\rho^\alpha (\alpha>1/2)$ was studied by Li, Li and Xin
 in \cite{LLX} and the phenomena of vacuum vanishing and
blow-up of solutions were found there. The global existence of
weak solutions for the initial-boundary-value problems for
spherically symmetric compressible Navier-Stokes equations with
density-dependent viscosity was proved by Guo, Jiu and Xin in
\cite{GJX}. More recently, there are some results on Cauchy
problem \eqref{(1.1)}-\eqref{(1.3)}. The existence and uniqueness
of global strong solutions to the compressible Navier-Stokes
equations \eqref{(1.1)}-\eqref{(1.3)} were obtained by Mellet and
Vasseur \cite{MV} where no vacuum is permitted in the initial
density. However, the a priori estimates obtained in \cite{MV}
depends on the time interval thus does not yield the
time-asymptotic behavior of the solutions. The first result about
the time-asymptotic behavior of the solutions to the Cauchy
problem \eqref{(1.1)}-\eqref{(1.3)} is obtained by Jiu-Xin
\cite{JX}, where the global existence, large time-asymptotic
behavior, the vanishing of the vacuum and the blow-up phenomena of
the weak solutions were considered in the case that $\r_+=\r_-$
and $u_+=u_-=0$. It is well-known that the large time-asymptotic
behavior of solutions to the system \eqref{(1.1)}-\eqref{(1.3)}
with different far field states of the initial values is closely
related to the corresponding Euler system \eqref{(1.4)}.
Rarefaction wave is one of the fundamental wave patterns to the
Euler system \eqref{(1.4)}. A natural question is how about the
stability of rarefaction waves to the compressible Navier-Stokes
equations \eqref{(1.1)}-\eqref{(1.3)} in the framework of weak
solutions which may contain vacuum states.

In the case $\r_\pm>0$, that is, the rarefaction wave does not
contain vacuum, we study in this paper the global existence, large
time-asymptotic behavior, vanishing of the vacuum and the blow-up
phenomena of weak solutions to the Cauchy problem
\eqref{(1.1)}-\eqref{(1.3)}. First we will construct a class of
approximate solutions satisfying some uniform estimates and
furthermore prove the global existence of weak solutions for
Cauchy problem of \eqref{(1.1)}-\eqref{(1.3)}. Moreover, the
time-asymptotic behaviors of weak solutions are investigated. More
precisely, it is proved that the density $\r$ tends to a
rarefaction wave as $t\to\infty$. As a consequence, there exists a
time $T_0>0$ such that when $t>T_0$, the vacuum states vanish and
the global weak solution becomes a unique strong one. Moreover,
the stability of the rarefaction wave is obtained in some weak
sense.

 \vskip 2mm
\noindent\emph{Notations.} Throughout this paper, positive generic
constants are denoted by $c$ and $C$, which are independent of
$\v$, $t$ and $T$, without confusion, and $C(\cdot)$ stands for
some generic constant(s) depending only on the quantity listed in
the parenthesis. For function spaces, $L^{p}(\Omega), 1\leq p\leq
\infty$, denote the usual Lebesgue spaces on $\Omega \subset
\mathbb{R}:=(-\infty,\infty)$. $W^{k,p}(\Omega)$ denotes the
$k^{th}$ order Sobolev space, $H^{k}(\Omega):=W^{k,2}(\Omega)$,
$\parallel \cdot
\parallel:=\parallel \cdot
\parallel_{L^{2}(\Omega)}$, and $\parallel \cdot
\parallel_k:=\parallel \cdot
\parallel_{H^{k}(\Omega)}$  for simplicity.  The domain
$\Omega$ will be often abbreviated without confusion.

\section{Preliminaries and Main Results}

\subsection{Rarefaction waves}

Consider   the solution to the following Cauchy problem for Burgers
equation
\begin{eqnarray}\label{(2.1)}
\begin{cases}
 \di w_{t}+ww_{x}=0,  \quad \,t>0,  x\in\mathbb{R}, \cr
 \di w_0( x):= w( 0,x)=\f{w_++w_-}{2}+\f{w_+-w_-}{2}K_q
\int^{\eta x }_0(1+y^2)^{-q}\,dy.
\end{cases}
\end{eqnarray}
Here $q\geq 2$ is some fixed constant, and $K_q$ is a constant such
that $K_q\int^{\infty}_0(1+y^2)^{-q} dy=1 $, and $\eta$ is a small
positive constant to be determined later. It is easy to see that the
solution to the above Burgers equation is given by
\begin{eqnarray}
w(t,x)=w_0(x_0(t, x)),\quad\quad x=x_0(t, x)+w_0(x_0(t, x))t.
\end{eqnarray}
Then the following properties hold (see
\cite{Matsumura-Nishihara-2}).

{\em Lemma 2.1}  Let $0\leq w_-<w_+$, Burgers equation
$\eqref{(2.1)}$ has a unique smooth solution $w(t,x)$ satisfying

i) ~ $w_-< w(t,x)<w_+,~w_x(t,x)> 0 $,

ii) ~ For any $p$ $(1\leq p\leq \infty)$, there exists a constant
$C_{pq}$ such that
\begin{eqnarray}\nonumber
&& \parallel w_x(t)\parallel_{L^p}\leq
C_{pq}\min\big{\{}\delta_r\eta^{1- \f1p},~
\delta_r^{\f1p}t^{-1+\f1p}\big{\}},   \cr&&\parallel
w_{xx}(t)\parallel_{L^p}\leq
C_{pq}\min\big{\{}\delta_r\eta^{2-\f1p},~
\eta^{(1-\f{1}{2q})(1-\f1p)}\d_r^{-\f{(p-1)}{2pq}}t^{-1-\f{(p-1)}{2pq}}\big{\}},
\end{eqnarray}
where $\d_r=w_+-w_-,$

 iii) ~ $\sup\limits_{x\in\mathbf{R}} |w(t,x)-w^r(\f
xt)|\rightarrow 0$, as $t\rightarrow \infty$.

Define a $2-$rarefaction wave $(\bar\r, \bar u) (t,x)$ by
\begin{equation}
\begin{array}{l}
\di \l_2(\r_\pm,u_\pm)=w_\pm,\quad w_-<w_+,\\
\di \l_2(\bar\r(t,x),\bar u(t,x))=w(1+t,x),\\
\di \Sigma_2(\bar\r(t,x),\bar u(t,x))=\Sigma_2(\r_{\pm},u_\pm),\\
\end{array}\label{(2.3)}
\end{equation}
Thus $(\bar\r, \bar u) (t,x)$ satisfies the system
\begin{equation}\label{(2.4)}
\left\{
\begin{array}{l}
\di \bar\r_t+(\bar\r\bar u)_x=0\\
\di (\bar\r \bar u)_t+(\bar\r\bar u^2+p(\bar\r))_x=0,\\
\end{array}
\right.
\end{equation}

 \vskip 2mm
\noindent{\em Lemma 2.2} \emph{ The $2-$rarefaction wave $(\bar\r,
\bar u) (t,x)$ satisfies}
\newcounter{w}
\begin{list}
{\upshape  \roman{w}) \;} {\setlength{\parsep}{\parskip}
 \setlength{\itemsep}{0ex plus0.1ex}
 \setlength{\rightmargin}{1em}
 \setlength{\leftmargin}{3em}
 \setlength{\labelwidth}{4em}
 \setlength{\labelsep}{0.1em}
 \usecounter{w}\setcounter{w}{0}
}
\item \emph{  $\bar\r_x>0,\quad \bar u_x>0$};
\item  \emph{For any $p$} ($1\leq p\leq \infty$), \emph{there exists a constant
$C_{pq}$ such that
\begin{eqnarray}\nonumber
&& \|(\bar\r_x, \bar u_x)(t,\cdot)\|_{L^p(\mathbf{R})}\leq
C_{pq}\min\{\d\eta^{1-\f1p}, \d^\f1p(1+t)^{-1+\f1p}\} ,
\cr&&\|(\bar\r_{xx},\bar u_{xx})(t,\cdot)\|_{L^p(\mathbf{R})}\leq
C_{pq}\min\{\d\eta^{2-\f1p},
\eta^{(1-\f{1}{2q})(1-\f1p)}\d^{-\f{p-1}{2pq}}(1+t)^{-1-\f{p-1}{2pq}}+\d^{\f1p}(1+t)^{-2+\f1p}\},
 \end{eqnarray}
 where $\d=|\r_+-\r_-|+|u_+-u_-|$ is the strength of the rarefaction
 wave;
 \item   $ \lim\limits_{t\rightarrow
\infty}\sup\limits_{\x\in\mathbf{R}}\big| (\bar\r,\bar
u)(t,x)-(\r^{r},u^{r})(\f x{1+t})\big| =0$.}
\end{list}

{\em Remark:} For any $1<p\leq +\i$,
$$
\int_0^T\|(\bar\r_{xx},\bar
u_{xx})(t,\cdot)\|_{L^p(\mathbf{R})}dt\leq C,
$$
where $C$ is independent of $T$. Note that in the case $p=1$, the
constant $C$ in the above estimates is not uniform in $T$.

Moreover, the following estimate holds:
$$
\int_0^T\|(\bar\r_{xx},\bar
u_{xx})(t,\cdot)\|_{L^\i(\mathbf{R})}dt\leq
C\eta^{\f2{4q+1}}\int_0^T(1+t)^{-1-\f1{4q+1}}dt\leq
C\eta^{\f2{4q+1}},
$$

\subsection{Main Results}
Set
\begin{equation}\label{10-21-3}
\begin{array}{ll}
\di \Psi(\r,\bar\r)&\di =\int_{\bar\r}^\r\f{p(s)-p(\bar\r)}{s^2}ds\\
&\di
=\f{1}{(\g-1)\r}\Big[\r^\g-\bar\r^\g-\g\bar\r^{\g-1}(\r-\bar\r)\Big].
\end{array}
\end{equation}
The initial data are assumed to satisfy:
 \begin{equation}
 \left\{
 \begin{array}{ll}
  &\rho_0\ge 0; \ \ \ m_0=0\ \ a.e. {\rm on}\ \{x\in \mathbf{R}|\rho_0(x)=0\};\\
 &
 \quad (\rho_0^{\alpha-\frac 12})_x\in
 L^2(\mathbf{R}),\quad \r_0\Psi(\r_0,\bar\r_0)\in L^1(\mathbf{R});\\
 &\di \r_0(\frac{m_0}{\rho_0}-\bar u_0)^2\in L^1(\mathbf{R}),\quad \r_0(\frac{m_0}{\r_0}-\bar u_0)^{3}\in
 L^1(\mathbf{R}),
 \end{array}
\right.\label{(2.5)}
 \end{equation}
where $\alpha>\frac 12$ and $(\bar\r_0,\bar u_0):=(\bar\r,\bar
u)(0,x)$ is the initial values of the $2-$rarefaction wave
$(\bar\r,\bar u)(t,x)$ constructed in section 2.1. Note that
\eqref{(2.5)} implies that $\rho_0\in C(\mathbf{R})$ which is the
space of continuous functions.

Before stating the main results, we give the definition of  weak
solutions to \eqref{(1.1)}-\eqref{(1.3)} associated with
$2-$rarefaction wave $(\bar\r,\bar u)(t,x)$ in \eqref{(2.4)}.

{\em Definition 2.1.} A pair $(\rho, u)$ is said to be a weak
solution to \eqref{(1.1)}-\eqref{(1.3)} towards the  rarefaction
wave $(\bar\r,\bar u)(t,x)$ in \eqref{(2.4)},  provided that

(1) $\rho\geq 0$ a.e., and
\begin{eqnarray*}
&&\rho\in L^\infty(0,T; L^\infty(\mathbf{R})))\cap C([0,\infty);W^{1,\infty}(\mathbf{R})^*), \\
&& (\rho^{\alpha-\frac12})_x\in L^\infty(0,T;L^2(\mathbf{R}))
,\sqrt{\rho}(u-\bar u)\in L^\infty(0,T;L^2(\mathbf{R})),
\end{eqnarray*}
where $W^{1,\infty}(\mathbf{R})^*$ is the dual space of
$W^{1,\infty}(\mathbf{R})$;

 (2) For any $t_2\geq
t_1\ge 0$ and any $\zeta\in C_0^1({\mathbf{R}}\times[t_1,t_2])$, the
mass equation \eqref{(1.1)} holds in the following sense:
\begin{eqnarray}
\int_{\mathbf{R}} (\rho-\bar\r)\zeta
dx|_{t_1}^{t_2}=\int_{t_1}^{t_2}\int_{\mathbf{R}}[(\rho-\bar\r)
\zeta_t+(\rho u-\bar\r\bar u)\cdot\zeta_x]dxdt; \label{mass}
\end{eqnarray}

(3) For any $\psi\in C_0^\infty({\mathbf{R}}\times[0,T))$, it holds
that
\begin{eqnarray}
&&\int_{\mathbf{R}} (m_0-\bar\r_0\bar u_0)\psi(0,\cdot)
dx+\int_0^T\int_{\mathbf{R}}\{\sqrt{\rho}[\sqrt{\rho}(u-\bar
u)]+(\r-\bar\r)\bar u\}\psi_t\\
&& +\{[\sqrt{\rho}(u-\bar u)]^2-2\sqrt{\r}\sqrt\r(u-\bar u)\bar
u+(\r-\bar \r)\bar u^2+(\rho^\gamma-\bar\r^\g)\}\psi_xdxdt \nonumber\\
&&+<\rho^\alpha (u-\bar u)_x,
\psi_x>+\int_0^T\int_{\mathbf{R}}\r^\a\bar u_x\psi_x dxdt=0,
\label{momentum}
\end{eqnarray}
where the diffusion term makes sense when written as
\begin{eqnarray}
&&<\rho^\alpha (u-\bar u)_x, \psi>=-\int_0^T\int_R
\rho^{\alpha-\frac12}\sqrt{\rho} (u-\bar u)\psi_{x}
dxdt\nonumber\\
&&-\frac{2\alpha}{2\alpha-1}\int_0^T\int_R
(\rho^{\alpha-\frac12})_x\sqrt{\rho} (u-\bar u)\psi dxdt.
\label{diffusion}
\end{eqnarray}

Our main results read as

{\em Theorem 2.1} (Existence of a weak solution) Let $\a$ and $\g$
satisfy that
\begin{equation}
\f12<\a\leq \f{\g+1}{2}.
\end{equation}
Suppose that \eqref{(2.5)} holds.
Then the Cauchy problem \eqref{(1.1)}-\eqref{(1.3)} admits a global
weak solution $(\rho(x,t), u(x,t))$ satisfying
\begin{eqnarray}
\rho\in  C(\mathbf{R}\times (0,T)),\label{th12-1}
\end{eqnarray}
\begin{eqnarray}
&&\rho\ge 0, \quad \max_{(x,t)\in
\mathbf{R}\times[0,T]}\rho\le C, \label{th12-2}\\
&&\sup_{t\in[0,T]}\int_\mathbf{R} (|\sqrt{\rho}(u-\bar
u)|^2+(\rho^{\alpha-\frac12})_x^2+\frac{1}{\gamma-1}[\rho^\gamma-\bar\rho^\gamma-\gamma\bar\rho^{\gamma-1}
(\rho-\bar\rho)]
dx\nonumber\\
&&+\int_0^T\int_\mathbf{R}
([(\rho^{\frac{\gamma+\alpha-1}{2}}-\bar\rho^{\frac{\gamma+\alpha-1}{2}})_x]^2+\Lambda(x,t)^2)
dxdt\le C,\label{th12-3}
\end{eqnarray}
where $C$ is an absolute constant depending on the initial data
and $\Lambda(x,t)\in L^2(\mb{R}\times (0,T))$ satisfying
\begin{eqnarray}
&&\int_0^T\int_\mathbf{R} \Lambda\varphi
dxdxt=-\int_0^T\int_\mathbf{R} \rho^{\alpha-\frac12}\sqrt{\rho}
(u-\bar u) \varphi_x
dxdt\nonumber\\
&&-\frac{2\alpha}{2\alpha-1}\int_0^T\int_\mathbf{R}
(\rho^{\alpha-\frac12})_x\sqrt{\rho} (u-\bar u)\varphi dxdt.
\label{diffusion1}
\end{eqnarray}

\noindent{\bf Remark 1:} Obviously, the case of shallow water
equation, i.e. $\a=1$, $\g=2$, is included in our theorem.

{\em Theorem 2.2} Suppose that $(\rho(x,t), u(x,t))$ is a weak
solution of the Cauchy problem \eqref{(1.1)}-\eqref{(1.3)}
satisfying \eqref{th12-1}, \eqref{th12-2} and \eqref{th12-3}. Then
we have
\begin{eqnarray}
\lim_{t\to\infty}\sup_{x\in \mathbf{R}}|\rho-\bar\rho|=0.
\label{th14-1}
\end{eqnarray}

Based on Theorem 2.2, it is easy to deduce that  under the
assumption that the approximate rarefaction wave satisfies
$\inf_{x,t}\bar\rho(t,x)>0$, there exists a time $T_0>0$ after which
the density has a positive lower bound and the vacuum states vanish.
Moreover, it will be shown that after the time $t=T_0$, the weak
solution becomes a unique strong one. Precisely, we have

{\em Theorem 2.3} Suppose that the assumptions of Theorem 2.1
hold.  Let $(\rho(x,t), u(x,t))$ be a weak solution of the Cauchy
problem \eqref{(1.1)}-\eqref{(1.3)} satisfying
\eqref{th12-1},\eqref{th12-2} and \eqref{th12-3}. Then for any
$0<\rho_1<\inf_{t,x}\bar\rho(t,x)$, there exists a time $T_0$ such
that
\begin{eqnarray}
0<\rho_1\le \rho(x,t)\le C, \ \ \ (x,t)\in \mathbf{R}\times
[T_0,\infty),
\end{eqnarray}
where $C$ is a constant same as in \eqref{th12-2}. Moreover, for
$t\ge T_0$, the weak solution becomes a unique strong solution to
\eqref{(1.1)}-\eqref{(1.3)}, satisfying
\begin{eqnarray*}
&&\rho-\bar\rho\in L^\infty(T_0, t;H^1(\mathbf{R})), \ \ \rho_t\in
L^\infty(T_0,t; L^2(\mathbf{R})),\\
&&u-\bar u\in L^2(T_0,t; H^2(\mathbf{R})), \ \ u_t\in L^2(T_0,t;
L^2(\mathbf{R}))
\end{eqnarray*}
and
\begin{eqnarray}
\sup_{x\in
\mathbf{R}}|\rho-\bar\rho|+\|\rho-\bar\rho\|_{L^p(\mathbf{R})}+\|u-\bar
u\|_{L^2(\mathbf{R})}\to 0 \label{Jan21-1}
\end{eqnarray}
as $t\to\infty$, where $2<p\le\infty$.

\noindent{\bf Remark 2:} It is interesting to note that there is
no requirement on the sizes of the strength of the rarefaction
wave and the perturbations. The class of initial perturbations
given by \eqref {(2.5)} is quite large compared with those for the
constant viscosity case, \cite{Liu-Xin-1},
\cite{Matsumura-Nishihara-1}, \cite{Matsumura-Nishihara-2}.

In addition, similar to \cite{LLX}, we can obtain some results on
the blow-up phenomena of the solution when the vacuum states
vanish, which will be presented in Section~5.

\section{Existence of a weak solution}

We first study the following approximate system:
\begin{eqnarray}\label{(3.1)}
\begin{cases}
    \rho   _{t}+ (\rho u) _{x} = 0,   \qquad\qquad\qquad\qquad\qquad\qquad x\in\mathbf{R},~~ t>0,\cr
    (\rho u)_{t}+  \big {(} \rho u^2 +p(\r)\big {)}_x
    =(\mu_\v(\r)  u_{x})_x,\cr
   \end{cases}
\end{eqnarray}
where $\mu_\v(\r)=\r^\a+\v\r^\t,\v>0,\t\in(0,\f12)$.

The initial values can be regularized in the following way. Due to
\eqref{(2.5)}, we have
$$
\lim_{x\rightarrow\pm\i}\r_0(x)=\r_{\pm}.
$$
For any suitably small constant $\l>0$, there exists $M>0$ such that
if $|x|\geq M$, then
$$
\r_0(x)\geq \l.
$$
For definiteness, we can choose and fix $\l=\f{\r_-}{2}$.

Define
$$
\r_0^1(x):=\left\{
\begin{array}{ll}
 \r_0(x)+\v^{\f{1}{2\a-2\t}},&\di {\rm if }~~|x|\leq M,\\
 \r_0(x)+\v^{\f{1}{2\a-2\t}}[(M+1)\mp x],\qquad &\di {\rm if}~~ M\leq\pm x\leq M+1,\\
 \r_0(x),&{\rm if}~~ |x|\geq M+1
\end{array}
\right.
$$
Then, $\r_0^1(x)$ is a continuous function in $\mathbf{R}$ and
$\r_0^1(x)\geq
\min\{\f{\r_-}{2},\v^{\f{1}{2\a-2\t}}\}=\v^{\f{1}{2\a-2\t}}$ for
suitably small $\v$. Also,
$$
0\leq\r_0^1(x)-\r_0(x)=(\r_0^1(x)-\r_0(x)){\bf 1}|_{|x|\leq M+1}\leq
2(M+1)\v^{\f{1}{2\a-2\t}}.
$$
Set
$$
\r_{0\v}(x)=(\r_0^1-\bar\r_0)\ast j_\v (x)+\bar\r_0(x).
$$
Hence, $\r_{0\v}\in C^\i(\mathbf{R})$ and $\di \lim_{x\rightarrow\pm
\i}\r_{0\v}=\r_{\pm}$ uniformly in $\v$. So for $\f{\r_-}{2}>0$,
there exists $M_1>0$ such that
$$
\r_{0\v}\geq \f{\r_-}{2},\qquad {\rm if}~~|x|\geq M_1.
$$
Note that $\r_{0\v}$ has a lower bound in the domain $|x|\leq
M_1.$ Indeed, since $(\r_0^1-\bar\r_0)\in C(\mathbf{R})$, it
follows that
$$
(\r_0^1-\bar\r_0)\ast j_\v\rightarrow  \r_0^1-\bar\r_0,\quad {\rm
in}~~C(|x|\leq M_1)
$$
uniformly as $\v\rightarrow 0.$ Hence
$$
\r_{0\v}\rightarrow \r_0^1(x)\quad {\rm in}~~C(|x|\leq M_1)
$$
uniformly as $\v\rightarrow 0.$

This implies that
$$
\r_{0\v}(x)\geq \f{1}{2}\v^{\f{1}{2\a-2\t}}\quad {\rm if} ~~|x|\leq
M_1.
$$
Thus
$$
\r_{0\v}(x)\geq
\min\{\f{\r_-}{2},\f{1}{2}\v^{\f1{2\a-2\t}}\}=\f12\v^{\f1{2\a-2\t}},
~~\forall x\in \mathbf{R}
$$
for suitably small $\v$.

Moreover,
$$
\r_{0\v}\Psi(\r_{0\v},\bar\r_0)\rightarrow
\r_{0}\Psi(\r_{0},\bar\r_0) ~{\rm in}~L^1(\mathbf{R}),\quad
(\rho_{0\v}^{\alpha-1/2})_x\to (\rho_{0}^{\alpha-1/2})_x \ {\rm in}\
L^2(\mathbf{R}).
$$
Therefore, it holds that
$$
\v^2[(\r_{0\v}^{\t-\f12})_x]^2=(\f{\t-\f12}{\a-\f12})^2\v^2\r_{0\v}^{2\t-2\a}[(\r_{0\v}^{\a-\f12})_x]^2\leq
C
$$
uniformly in $\v$.

To regularize $m_0$, one sets
$$
F_0=\r_0(\f{m_0}{\r_0}-\bar u_0)^{3}\in L^1(\mathbf{R})
$$
and
$$
F_{0\v}=F_0\ast j_\v.
$$
Then,
$$F_{0\v}\rightarrow F_0, ~~~{\rm a.e.}$$
and
$$F_{0\v}\rightarrow F_0,~~ {\rm in}~~ L^1(\mathbf{R}).$$
Now we define
$$
m_{0\v}=\r_{0\v}\Big(\bar
u_0+\left(\f{F_{0\v}}{\r_{0\v}}\right)^{\f{1}{3}}\Big).
$$
Then
$$
\r_{0\v}(\f{m_{0\v}}{\r_{0\v}}-\bar
u_0)^{3}\rightarrow\r_0(\f{m_0}{\r_0}-\bar u_0)^{3}~~{\rm in}~~
L^1(\mathbf{R})
$$
and
$$
\r_{0\v}(\f{m_{0\v}}{\r_{0\v}}-\bar
u_0)^{2}\rightarrow\r_0(\f{m_0}{\r_0}-\bar u_0)^{2}~~{\rm in}~~
L^1(\mathbf{R}).
$$

For any fixed $T>0$ and for any fixed $\v>0$, there exists a
unique smooth approximate solution to (\ref{(3.1)}) in the region
$(x,t)\in \mathbf{R}\times(0,T)$ with initial data
\begin{equation}
(\r,\r u)(0,x)=(\r_{0\v},m_{0\v}).\label{10-21-1}
\end{equation}
We refer to \cite{MV2} for the wellposedness of the global strong
solution to the approximate system (\ref{(3.1)}) with
\eqref{10-21-1}.

 The following estimates are crucial to
prove our main results, which are based on the usual energy
estimates and a new entropy estimates (see \cite{BD1}-\cite{BD3}).

{\em Lemma 3.1} Let
\begin{equation}
\f12<\a\leq \f{\g+1}{2},\label{(3.2)}
\end{equation}
Suppose that $(\rho_\v, u_\v)$ is a smooth solution to
\eqref{(3.1)} satisfying $\rho_\epsilon>0$. Then for any $T>0$ and
$\epsilon>0$ satisfying $\sqrt\v\ln(1+T)\leq \v^{\f14}$, the
following estimates hold:
\begin{equation}
\begin{array}{ll}
\di \sup_{t\in[0,T]}\int_{\mathbf{R}}\Big\{\r_\v(u_\v-\bar
u)^2+\Big[\Big(\f{\r_\v^{\a-\f12}}{\a-\f12}\Big)_x\Big]^2+\v^2\Big[\Big(\f{\r_\v^{\t-\f12}}{\t-\f12}\Big)_x\Big]^2+\r_\v\Psi(\r_\v,\bar\r)\Big\}(x,t)dx\\
\di +\int_0^T\int_{\mathbf{R}}\Big\{\bar
u_x\Big[p(\r_\v)-p(\bar\r)-p^\prime(\bar\r)(\r_\v-\bar\r)\Big]+\r_\v(u_\v-\bar
u)^2\bar u_x+(\r_\v^\a+\v\r_\v^\t)\Big[(u_\v-\bar u)_x\Big]^2\\
\di +\Big[(\r_\v^{\f{\a+\g-1}{2}}-\bar\r^{\f{\a+\g-1}{2}})_x\Big]^2
+\v\Big[(\r_\v^{\f{\t+\g-1}{2}}-\bar\r^{\f{\t+\g-1}{2}})_x\Big]^2\Big\}(x,t)dxdt\\
\di \leq C,
\end{array}\label{(3.3)}
\end{equation}
where $C$ is an universal constant independent of $\v$ and $T$.

In the following, the subscript $\v$ in the approximate solution
$(\r_\v,u_\v)(t,x)$ will be omitted for simplicity.

 {\em
Proof:} {\it Step 1.  Energy Equality}

 It follows from $(\ref{(3.1)})_2$ that
\begin{equation}
\r u_t+\r uu_x+p(\r)_x=(\mu_\v(\r)u_x)_x.\label{10-21-2}
\end{equation}
Subtracting \eqref{10-21-2} from the second equation of
\eqref{(2.4)} gives
\begin{equation}
\r(u-\bar u)_t+\r u(u-\bar u)_x+(p(\r)-p(\bar\r))_x+(\r-\bar\r)\bar
u_t+(\r u-\bar\r\bar u)\bar u_x=(\mu_\v(\r)(u-\bar
u)_x)_x+(\mu_\v(\r)\bar u_x)_x,\label{(3.5)}
\end{equation}
Multiplying \eqref{(3.5)} by $u-\bar u$, we  get
\begin{equation}
\begin{array}{ll}
\di \Big[\f{\r(u-\bar u)^2}{2}\Big]_t+\Big[\f{\r u(u-\bar
u)^2}{2}\Big]_x+(u-\bar u)(p(\r)-p(\bar\r))_x-\Big[\mu_\v(\r)(u-\bar
u)(u-\bar u)_x\Big]_x\\[2mm]
\di +\mu_\v(\r)\Big[(u-\bar u)_x\Big]^2=\Big[\mu_\v(\r)\Big]_x\bar
u_x(u-\bar u)+\mu_\v(\r)\bar u_{xx}(u-\bar u)-\Big[(\r-\bar\r)\bar
u_t+(\r u-\bar\r\bar u)\bar u_x\Big](u-\bar u).
\end{array}\label{(3.6)}
\end{equation}
Note that $\Psi(\r,\bar\r)$ defined in \eqref{10-21-3} satisfies
\begin{equation}
\begin{array}{l}
\di \Big[\r\Psi(\r,\bar\r)\Big]_t+\Big[\r
u\Psi(\r,\bar\r)\Big]_x+(u-\bar u)_x(p(\r)-p(\bar\r))+\bar
u_x\Big[\r^\g-\bar\r^\g-\g\bar\r^{\g-1}(\r-\bar\r)\Big]\\
\di =-\f{p(\bar\r)_x}{\bar\r}(\r-\bar\r)(u-\bar u).
\end{array}\label{(3.7)}
\end{equation}
It follows from \eqref{(3.6)} and \eqref{(3.7)} that
\begin{equation}
\begin{array}{ll}
\di \Big[\f{\r(u-\bar u)^2}{2}+\r\Psi(\r,\bar\r)\Big]_t+H_{1x}(t,x)
+\mu_\v(\r)\Big[(u-\bar u)_x\Big]^2+\bar
u_x\Big[\r^\g-\bar\r^\g-\g\bar\r^{\g-1}(\r-\bar\r)\Big]\\\di
=\Big[\mu_\v(\r)\Big]_x\bar u_x(u-\bar u)+\mu_\v(\r)\bar
u_{xx}(u-\bar u)-\Big[(\r-\bar\r)\bar u_t+(\r u-\bar\r\bar u)\bar
u_x+\f{p(\bar\r)_x}{\bar\r}(\r-\bar\r)\Big](u-\bar u),
\end{array}\label{(3.8)}
\end{equation}
where
$$
H_1(t,x)=\f{\r u(u-\bar u)^2}{2}+\r u\Psi(\r,\bar\r)+(u-\bar
u)(p(\r)-p(\bar\r))-\mu_\v(\r)(u-\bar u)(u-\bar u)_x.
$$
Since
$$
(\r-\bar\r)\bar u_t+(\r u-\bar\r\bar u)\bar
u_x+\f{p(\bar\r)_x}{\bar\r}(\r-\bar\r)=\r(u-\bar u)\bar u_x,
$$
we obtain
\begin{equation}
\begin{array}{ll}
\di \Big[\f{\r(u-\bar u)^2}{2}+\r\Psi(\r,\bar\r)\Big]_t+H_{1x}(t,x)
+\mu_\v(\r)\Big[(u-\bar u)_x\Big]^2+\bar
u_x\Big[\r^\g-\bar\r^\g-\g\bar\r^{\g-1}(\r-\bar\r)\Big]\\\di
+\r(u-\bar u)^2\bar u_x=\Big[\mu_\v(\r)\Big]_x\bar u_x(u-\bar
u)+\mu_\v(\r)\bar u_{xx}(u-\bar u).
\end{array}\label{(3.9)}
\end{equation}

{\it Step 2. New Entropy Equality}

Rewrite \eqref{(3.5)} as
\begin{equation}
\r(u-\bar u)_t+\r u(u-\bar u)_x+(p(\r)-p(\bar\r))_x+(\r-\bar\r)\bar
u_t+(\r u-\bar\r\bar u)\bar u_x=\Big[(\r^{\a-1}+\v\r^{\t-1})\r
u_x\Big]_x.\label{(3.10)}
\end{equation}
Note that
\begin{equation}
\di \Big[(\r^{\a-1}+\v\r^{\t-1})\r u_x\Big]_x \di
=-\r(\varphi_\v^{\a,\t}(\r))_{xt}-\r
u(\varphi_\v^{\a,\t}(\r))_{xx},\label{(3.11)}
\end{equation}
where $\varphi_\v^{\a,\t}(\r),~0<\t<\f12,$ is defined by
$$
\varphi_\v^{\a,\t}(\r)=\left\{
\begin{array}{ll}
\di \f{\r^{\a-1}}{\a-1}+\v\f{\r^{\t-1}}{\t-1},&\di {\rm
if}~~\a\neq1,\a>0,\\
\di \ln\r+\v\f{\r^{\t-1}}{\t-1},&\di {\rm if}~~\a=1.
\end{array}
\right.
$$
Thus \eqref{(3.10)} becomes
\begin{equation}
\r(u-\bar u)_t+\r u(u-\bar u)_x+(p(\r)-p(\bar\r))_x+(\r-\bar\r)\bar
u_t+(\r u-\bar\r\bar u)\bar u_x=-\r(\varphi_\v^{\a,\t}(\r))_{xt}-\r
u(\varphi_\v^{\a,\v}(\r))_{xx}.\label{(3.12)}
\end{equation}
Multiplying \eqref{(3.12)} by $(\varphi_\v^{\a,\t}(\r))_x$ shows
that
\begin{equation}
\begin{array}{ll}
\di \Big[\f{\r(\varphi_\v^{\a,\t}(\r))_x^2}{2}\Big]_t+ \Big[\f{\r
u(\varphi_\v^{\a,\t}(\r))_x^2}{2}\Big]_x+\Big[\r(u-\bar
u)(\varphi_\v^{\a,\t}(\r))_x\Big]_t +\Big[\r u(u-\bar
u)(\varphi_\v^{\a,\t}(\r))_x\Big]_x\\
\di -(u-\bar u)\Big[\r(\varphi_\v^{\a,\t}(\r))_{xt}+\r
u(\varphi_\v^{\a,\t}(\r))_{xx}\Big]+(\varphi_\v^{\a,\t}(\r))_x(p(\r)-p(\bar\r))_x\\
\di +(\varphi_\v^{\a,\t}(\r))_x\Big[(\r-\bar\r)\bar u_t+(\r
u-\bar\r\bar u)\bar u_x\Big]=0.
\end{array}\label{(3.13)}
\end{equation}
Combining \eqref{(3.12)} with \eqref{(3.13)} yields
\begin{equation}
\begin{array}{ll}
\di \Big\{\f12\r\left[(u-\bar
u)+(\varphi_\v^{\a,\t}(\r))_x\right]^2\Big\}_t+ \Big\{\f12\r
u\left[(u-\bar
u)+(\varphi_\v^{\a,\t}(\r))_x\right]^2\Big\}_x+(u-\bar u)(p(\r)-p(\bar\r))_x\\
\di +(\varphi_\v^{\a,\t}(\r))_x(p(\r)-p(\bar\r))_x+(u-\bar
u)\Big[(\r-\bar\r)\bar u_t+(\r
u-\bar\r\bar u)\bar u_x\Big]\\
\di +(\varphi_\v^{\a,\t}(\r))_x\Big[(\r-\bar\r)\bar u_t+(\r
u-\bar\r\bar u)\bar u_x\Big]=0.
\end{array}\label{(3.14)}
\end{equation}

{\it Step 3. A Priori Estimates}

It follows from \eqref{(3.7)} and \eqref{(3.14)} that
\begin{equation}
\begin{array}{ll}
\di \Big\{\f12\r\left[(u-\bar
u)+(\varphi_\v^{\a,\t}(\r))_x\right]^2+\r\Psi(\r,\bar\r)\Big\}_t+
\Big\{\f12\r u\left[(u-\bar
u)+(\varphi_\v^{\a,\t}(\r))_x\right]^2+\r u\Psi(\r,\bar\r)\\
\di +(u-\bar u)(p(\r)-p(\bar\r))\Big\}_x+\bar
u_x\Big[p(\r)-p(\bar\r)-p^\prime(\bar\r)(\r-\bar\r)\Big]+\r(u-\bar u)^2\bar u_x\\
\di +(\varphi_\v^{\a,\t}(\r))_x\Big[(\r-\bar\r)\bar u_t+(\r
u-\bar\r\bar u)\bar u_x+p(\r)_x-p(\bar\r)_x\Big]=0.
\end{array}\label{(3.15)}
\end{equation}
Now we deal with the last term on the left hand side of
\eqref{(3.15)}. Note that
\begin{equation}
\di (\r-\bar\r)\bar u_t+(\r u-\bar\r\bar u)\bar
u_x+p(\r)_x-p(\bar\r)_x \di =\r(u-\bar u)\bar
u_x+\Big[p(\r)_x-\f{\r p(\bar\r)_x}{\bar\r}\Big],\label{(3.16)}
\end{equation}
and
\begin{equation}
(\varphi_\v^{\a,\t}(\r))_x=\r^{\a-2}\r_x+\v\r^{\t-2}\r_x.\label{(3.17)}
\end{equation}
Thus
\begin{equation}
\begin{array}{ll}
&\di (\varphi_\v^{\a,\t}(\r))_x\Big[(\r-\bar\r)\bar u_t+(\r
u-\bar\r\bar u)\bar u_x+p(\r)_x-p(\bar\r)_x\Big]\\
&\di =\Big(\f{\r^\a}{\a}+\v\f{\r^\t}{\t}\Big)_x(u-\bar u)\bar
u_x+(\r^{\a-2}\r_x+\v\r^{\t-2}\r_x)\Big[p(\r)_x-\f{\r
p(\bar\r)_x}{\bar\r}\Big].\label{(3.18)}
\end{array}
\end{equation}
Direct computations show
\begin{equation}
\begin{array}{ll}
&\di \r^{\a-2}\r_x\Big[p(\r)_x-\f{\r p(\bar\r)_x}{\bar\r}\Big]\\
&\di
=\f{4\g}{(\a+\g-1)^2}\Big[(\r^{\f{\a+\g-1}{2}}-\bar\r^{\f{\a+\g-1}{2}})_x\Big]^2+\Big[\f{8\g}{(\a+\g-1)^2}(\bar\r^{\f{\a+\g-1}{2}})_x(\r^{\f{\a+\g-1}{2}}-\bar\r^{\f{\a+\g-1}{2}})\\
&\di \quad-\f{2\g}{\a(\a+\g-1)}(\bar\r^{\f{\a+\g-1}{2}})_x\bar\r^{\f{\g-\a-1}{2}}(\r^\a-\bar\r^\a)\Big]_x-\f{8\g}{(\a+\g-1)^2}(\bar\r^{\f{\a+\g-1}{2}})_{xx}(\r^{\f{\a+\g-1}{2}}-\bar\r^{\f{\a+\g-1}{2}})\\
&\di \quad+\f{2\g}{\a(\a+\g-1)}\Big[(\bar\r^{\f{\a+\g-1}{2}})_x\bar\r^{\f{\g-\a-1}{2}}\Big]_x(\r^\a-\bar\r^\a).\\
\end{array}\label{(3.19)}
\end{equation}
and
\begin{equation}
\begin{array}{ll}
&\di \r^{\t-2}\r_x\Big[p(\r)_x-\f{\r p(\bar\r)_x}{\bar\r}\Big]\\
&\di
=\f{4\g}{(\t+\g-1)^2}\Big[(\r^{\f{\t+\g-1}{2}}-\bar\r^{\f{\t+\g-1}{2}})_x\Big]^2+\Big[\f{8\g}{(\t+\g-1)^2}(\bar\r^{\f{\t+\g-1}{2}})_x(\r^{\f{\t+\g-1}{2}}-\bar\r^{\f{\t+\g-1}{2}})\\
&\di \quad-\f{2\g}{\t(\t+\g-1)}(\bar\r^{\f{\t+\g-1}{2}})_x\bar\r^{\f{\g-\t-1}{2}}(\r^\t-\bar\r^\t)\Big]_x-\f{8\g}{(\t+\g-1)^2}(\bar\r^{\f{\t+\g-1}{2}})_{xx}(\r^{\f{\t+\g-1}{2}}-\bar\r^{\f{\t+\g-1}{2}})\\
&\di \quad+\f{2\g}{\t(\t+\g-1)}\Big[(\bar\r^{\f{\t+\g-1}{2}})_x\bar\r^{\f{\g-\t-1}{2}}\Big]_x(\r^\t-\bar\r^\t).\\
\end{array}\label{(3.20)}
\end{equation}
Substituting \eqref{(3.18)}-\eqref{(3.20)} into \eqref{(3.15)}
gives
\begin{equation}
\begin{array}{ll}
\di \Big\{\f12\r\left[(u-\bar
u)+(\varphi_\v^{\a,\t}(\r))_x\right]^2+\r\Psi(\r,\bar\r)\Big\}_t+
H_{2x}(t,x)+\bar
u_x\Big[p(\r)-p(\bar\r)-p^\prime(\bar\r)(\r-\bar\r)\Big]\\
\di +\r(u-\bar u)^2\bar
u_x+\Big(\f{\r^\a}{\a}+\v\f{\r^\t}{\t}\Big)_x(u-\bar u)\bar
u_x+\f{4\g}{(\a+\g-1)^2}\Big[(\r^{\f{\a+\g-1}{2}}-\bar\r^{\f{\a+\g-1}{2}})_x\Big]^2\\
\di
+\v\f{4\g}{(\t+\g-1)^2}\Big[(\r^{\f{\t+\g-1}{2}}-\bar\r^{\f{\t+\g-1}{2}})_x\Big]^2=\f{8\g}{(\a+\g-1)^2}(\bar\r^{\f{\a+\g-1}{2}})_{xx}(\r^{\f{\a+\g-1}{2}}-\bar\r^{\f{\a+\g-1}{2}})\\
\di
+\v\f{8\g}{(\t+\g-1)^2}(\bar\r^{\f{\t+\g-1}{2}})_{xx}(\r^{\f{\t+\g-1}{2}}-\bar\r^{\f{\t+\g-1}{2}})-\f{2\g}{\a(\a+\g-1)}\Big[(\bar\r^{\f{\a+\g-1}{2}})_x\bar\r^{\f{\g-\a-1}{2}}\Big]_x(\r^\a-\bar\r^\a)\\
\di
-\v\f{2\g}{\t(\t+\g-1)}\Big[(\bar\r^{\f{\t+\g-1}{2}})_x\bar\r^{\f{\g-\t-1}{2}}\Big]_x(\r^\t-\bar\r^\t),
\end{array}\label{(3.21)}
\end{equation}
where
\begin{equation}
\begin{array}{ll}
\di H_2(t,x)=&\di\f12\r u\left[(u-\bar
u)+(\varphi_\v^{\a,\t}(\r))_x\right]^2+\r u\Psi(\r,\bar\r)+(u-\bar
u)(p(\r)-p(\bar\r))\\
&\di
+\f{8\g}{(\a+\g-1)^2}(\bar\r^{\f{\a+\g-1}{2}})_x(\r^{\f{\a+\g-1}{2}}-\bar\r^{\f{\a+\g-1}{2}})\\
&\di
-\f{2\g}{\a(\a+\g-1)}(\bar\r^{\f{\a+\g-1}{2}})_x\bar\r^{\f{\g-\a-1}{2}}(\r^\a-\bar\r^\a)\\
&\di +\v\f{8\g}{(\t+\g-1)^2}(\bar\r^{\f{\t+\g-1}{2}})_x(\r^{\f{\t+\g-1}{2}}-\bar\r^{\f{\t+\g-1}{2}})\\
&\di
-\v\f{2\g}{\t(\t+\g-1)}(\bar\r^{\f{\t+\g-1}{2}})_x\bar\r^{\f{\g-\t-1}{2}}(\r^\t-\bar\r^\t).
\end{array}\label{(3.22)}
\end{equation}
 Multiplying \eqref{(3.21)} by $\a$ and  then
adding up to \eqref{(3.9)} and noticing that
$\Big[\mu_\v(\r)\Big]_x=(\r^\a)_x+\v(\r^\t)_x$ in the right hand
side of \eqref{(3.9)}, we get
\begin{equation}
\begin{array}{ll}
\di \Big\{\f{\a}{2}\r\left[(u-\bar
u)+(\varphi_\v^{\a,\t}(\r))_x\right]^2+\f{\r(u-\bar
u)^2}{2}+(\a+1)\r\Psi(\r,\bar\r)\Big\}_t+ \Big[\a
H_{2}(t,x)+H_1(t,x)\Big]_x\\
\di +(\a+1)\bar
u_x\Big[p(\r)-p(\bar\r)-p^\prime(\bar\r)(\r-\bar\r)\Big]+(\a+1)\r(u-\bar
u)^2\bar u_x+(\r^\a+\v\r^\t)\Big[(u-\bar u)_x\Big]^2\\
\di
+\f{4\a\g}{(\a+\g-1)^2}\Big[(\r^{\f{\a+\g-1}{2}}-\bar\r^{\f{\a+\g-1}{2}})_x\Big]^2
+\v\f{4\a\g}{(\t+\g-1)^2}\Big[(\r^{\f{\t+\g-1}{2}}-\bar\r^{\f{\t+\g-1}{2}})_x\Big]^2\\
\di=\r^\a\bar u_{xx}(u-\bar u)+\v\Big[\r^\t\bar u_{xx}(u-\bar u)+(1-\f\a\t)(\r^\t)_x(u-\bar u)\bar u_x\Big]\\
\di +\f{8\a\g}{(\a+\g-1)^2}(\bar\r^{\f{\a+\g-1}{2}})_{xx}(\r^{\f{\a+\g-1}{2}}-\bar\r^{\f{\a+\g-1}{2}})\\
\di
+\v\f{8\a\g}{(\t+\g-1)^2}(\bar\r^{\f{\t+\g-1}{2}})_{xx}(\r^{\f{\t+\g-1}{2}}-\bar\r^{\f{\t+\g-1}{2}})
-\f{2\g}{(\a+\g-1)}\Big[(\bar\r^{\f{\a+\g-1}{2}})_x\bar\r^{\f{\g-\a-1}{2}}\Big]_x(\r^\a-\bar\r^\a)\\
\di
-\v\f{2\a\g}{\t(\t+\g-1)}\Big[(\bar\r^{\f{\t+\g-1}{2}})_x\bar\r^{\f{\g-\t-1}{2}}\Big]_x(\r^\t-\bar\r^\t)\\
\di := \sum_{i=1}^6I_i.
\end{array}\label{(3.22)}
\end{equation}
Integrating \eqref{(3.22)} over $[0,t]\times\mathbf{R}$ with respect
to $t,x$ gives
\begin{equation}
\begin{array}{ll}
\di \int_{\mathbf{R}}\Big\{\f{\a}{2}\r\left[(u-\bar
u)+(\varphi_\v^{\a,\t}(\r))_x\right]^2+\f{\r(u-\bar
u)^2}{2}+(\a+1)\r\Psi(\r,\bar\r)\Big\}(t,x)dx \\
 \di +\int_0^t\int_{\mathbf{R}}\Big\{(\a+1)\bar
u_x\Big[p(\r)-p(\bar\r)-p^\prime(\bar\r)(\r-\bar\r)\Big]+(\a+1)\r(u-\bar
u)^2\bar u_x\\
\di+(\r^\a+\v\r^\t)\Big[(u-\bar u)_x\Big]^2
+\f{4\a\g}{(\a+\g-1)^2}\Big[(\r^{\f{\a+\g-1}{2}}-\bar\r^{\f{\a+\g-1}{2}})_x\Big]^2\\
\di
+\v\f{4\a\g}{(\t+\g-1)^2}\Big[(\r^{\f{\t+\g-1}{2}}-\bar\r^{\f{\t+\g-1}{2}})_x\Big]^2\Big\}dxd\tau\\
\di = \int_0^t\int_{\mathbf{R}}\sum_{i=1}^6 I_idxd\tau.
\end{array}\label{(3.23)}
\end{equation}

We now estimate the right hand side of \eqref{(3.23)} terms by
terms. First,
\begin{equation}
\begin{array}{ll}
\di \int_0^t\int_{\mathbf{R}}I_1dxd\tau&\di
=\int_0^t\int_{\mathbf{R}}\r^\a\bar
u_{xx}(u-\bar u)dxd\tau\\
&\di =\int_0^t\int_{\mathbf{R}}\sqrt{\r}(u-\bar u)\r^{\a-\f12}\bar
u_{xx}dxd\tau\\
&\di =\int_0^t\int_{\mathbf{R}}\sqrt{\r}(u-\bar u)\r^{\a-\f12}\bar
u_{xx}[{\bf 1}|_{\{0\leq\r\leq2\r_+\}}+{\bf
1}|_{\{\r\geq2\r_+\}}]dxd\tau\\
&\di :=J_1^1+J_1^2,
\end{array}\label{(3.24)}
\end{equation}
where ${\bf 1}|_\Omega$ is the characteristic function of a set
$\Omega\subset(0,t)\times\mathbf{R}$.

Using Lemma 2.2 (and its Remark), noting that $\a>\f12$,  we have
\begin{equation}
\begin{array}{ll}
\di J_1^1&\di \leq C\int_0^t\|\sqrt\r(u-\bar
u)\|_{L^2(\mathbf{R})}\|\bar
u_{xx}\|_{L^2(\mathbf{R})}d\tau\\
&\di\leq C\sup_{t\in[0,T]}\|\sqrt\r(u-\bar
u)\|_{L^2(\mathbf{R})}\int_0^t\|\bar
u_{xx}\|_{L^2(\mathbf{R})}d\tau\\
&\di \leq C\sup_{t\in[0,T]}\|\sqrt\r(u-\bar u)\|_{L^2(\mathbf{R})},
\end{array}\label{(3.25)}
\end{equation}
and
\begin{equation}
\begin{array}{ll}
\di J_1^2&\di= \int_0^t\int_{\mathbf{R}}\sqrt{\r}(u-\bar
u)(\r^{\a-\f12}-\bar\r^{\a-\f12})\bar u_{xx}{\bf
1}|_{\{\r\geq2\r_+\}}dxd\tau+\int_0^t\int_{\mathbf{R}}\sqrt{\r}(u-\bar
u)\bar\r^{\a-\f12}\bar u_{xx}{\bf 1}|_{\{\r\geq2\r_+\}}dxd\tau\\
&\di\leq C\sup_{t\in[0,T]}\|\sqrt\r(u-\bar
u)\|_{L^2(\mathbf{R})}\sup_{t\in[0,T]}\|(\r^{\a-\f12}-\bar\r^{\a-\f12}){\bf
1}|_{\{\r\geq2\r_+\}}\|_{L^2(\mathbf{R})}\int_0^t\|\bar
u_{xx}\|_{L^\i(\mathbf{R})}d\tau\\
&\di \quad+C\sup_{t\in[0,T]}\|\sqrt\r(u-\bar
u)\|_{L^2(\mathbf{R})}\int_0^t\|\bar
u_{xx}\|_{L^2(\mathbf{R})}d\tau\\
&\di \leq C\eta^{\f{2}{4q+1}}\sup_{t\in[0,T]}\|\sqrt\r(u-\bar
u)\|_{L^2(\mathbf{R})}\sup_{t\in[0,T]}\|(\r^{\a-\f12}-\bar\r^{\a-\f12}){\bf
1}|_{\{\r\geq2\r_+\}}\|_{L^2(\mathbf{R})}\int_0^t(1+\tau)^{-1-\f{1}{4q+1}}d\tau\\
&\di \quad+ C\sup_{t\in[0,T]}\|\sqrt\r(u-\bar
u)\|_{L^2(\mathbf{R})}\\
&\di \leq C\eta^{\f{2}{4q+1}}\Big[\sup_{t\in[0,T]}\|\sqrt\r(u-\bar
u)\|^2_{L^2(\mathbf{R})}+\sup_{t\in[0,T]}\|(\r^{\a-\f12}-\bar\r^{\a-\f12}){\bf
1}|_{\{\r\geq2\r_+\}}\|^2_{L^2(\mathbf{R})}\Big]+C_{\eta}.
\end{array}\label{(3.26)}
\end{equation}

Note that if $\a$ and $\g$ satisfy
\begin{equation}
0<2(\a-\f12)\leq\g,\quad{\rm i.e.,}\quad \f12<\a\leq
\f{\g+1}{2},\label{(3.29)}
\end{equation}
then
\begin{equation}
\begin{array}{ll}
\di
\lim_{\r\rightarrow+\i}\f{(\r^{\a-\f12}-\bar\r^{\a-\f12})^2}{\r\Psi(\r,\bar\r)}&\di
=\lim_{\r\rightarrow+\i}\f{(\g-1)(\r^{\a-\f12}-\bar\r^{\a-\f12})^2}{\r^\g-\bar\r^\g-\g\bar\r^{\g-1}(\r-\bar\r)}\\
&\leq C.
\end{array}\label{(3.30)}
\end{equation}

Thus if $\f12<\a\leq \f{\g+1}{2}$, then
\begin{equation}
\sup_{t\in[0,T]}\|(\r^{\a-\f12}-\bar\r^{\a-\f12}){\bf
1}|_{\{\r\geq2\r_+\}}\|^2_{L^2(\mathbf{R})}\leq
C\sup_{t\in[0,T]}\|\r\Psi(\r,\bar\r)\|_{L^1(\mathbf{R})}\label{(3.31)}
\end{equation}
for some uniform constant $C>0$.

Combining \eqref{(3.25)}-\eqref{(3.31)} together shows that
\begin{equation}
\int_0^t\int_{\mathbf{R}}I_1dxd\tau\leq
C\eta^{\f{2}{4q+1}}\Big[\sup_{t\in[0,T]}\|\sqrt\r(u-\bar
u)\|^2_{L^2(\mathbf{R})}+\sup_{t\in[0,T]}\|\r\Psi(\r,\bar\r)\|_{L^1(\mathbf{R})}\Big]+C_{\eta}.
\end{equation}
Next, we estimate $\int_0^t\int_{\mathbf{R}}I_2dxd\tau$ which can
be rewritten as
$$
\begin{array}{ll}
\di \int_0^t\int_{\mathbf{R}}I_2dxd\tau&\di
=\v\int_0^t\int_{\mathbf{R}}\Big[\r^\t\bar
u_{xx}(u-\bar u)+(1-\f\a\t)(\r^\t)_x(u-\bar u)\bar u_x\Big]dxd\tau\\
&\di =-\v\int_0^t\int_{\mathbf{R}}\Big[\r^\t\bar u_{x}(u-\bar
u)_x+\f\a\t(\r^\t)_x(u-\bar u)\bar u_x\Big]dxd\tau\\
&\di :=J_2^1+J_2^2,
\end{array}
$$
Using Young inequality, one has
$$
\begin{array}{ll}
\di J_2^1&\di =-\v\int_0^t\int_{\mathbf{R}}\r^\t\bar u_{x}(u-\bar
u)_xdxd\tau\\
&\di \leq \f{\v}{4}\int_0^t\int_{\mathbf{R}} \r^\t[(u-\bar
u)_x]^2dxd\tau+ \v \int_0^t\int_{\mathbf{R}} \r^\t\bar u_x^2dxd\tau.
\end{array}
$$
By Lemma 2.2, one can obtain
$$
\begin{array}{ll}
\di \v \int_0^t\int_{\mathbf{R}} \r^\t\bar u_x^2dxd\tau&\di =\v
\int_0^t\int_{\mathbf{R}} \r^\t({\bf 1}_{\{0\leq \r\leq
2\r_+\}}+{\bf 1}_{\{\r\geq 2\r_+\}})\bar
u_x^2dxd\tau\\
&\di \leq C\v \int_0^t\int_{\mathbf{R}} \bar u_x^2dxd\tau+\v
\int_0^t\int_{\mathbf{R}} (\r^\t-\bar\r^\t){\bf 1}_{\{\r\geq
2\r_+\}}\bar u_x^2dxd\tau\\
&\di \leq C\v\ln(1+T)+C\v \sup_{t\in [0,T]}\|(\r^\t-\bar\r^\t){\bf
1}_{\{\r\geq 2\r_+\}}\|_{L^1(\mathbf{R})}\int_0^t\|\bar
u_x\|^2_{L^{\i}(\mathbf{R})}d\tau\\
&\di \leq C\v\ln(1+T)+C\v \sup_{t\in
[0,T]}\|\r\Psi(\r,\bar\r)\|_{L^1(\mathbf{R})}
\end{array}
$$
due to the fact that
$$
\lim_{\r\ra+\i}\f{|(\r^\t-\bar\r^\t)|{\bf 1}_{\{\r\geq
2\r_+\}}}{\r\Psi(\r,\bar\r)}=0.
$$

Moreover, direct estimates show
$$
\begin{array}{ll}
\di J_2^2&\di
=-\f{\v\t}{\t-\f12}\int_0^t\int_{\mathbf{R}}(\r^{\t-\f12})_x\bar
u_{x}\sqrt\r(u-\bar
u)dxd\tau\\
&\di \leq
C\sqrt\v\sup_{t\in[0,T]}\|\sqrt\v(\r^{\t-\f12})_x\|_{L^2(\mathbf{R})}\sup_{t\in[0,T]}\|\sqrt\r(u-\bar
u)\|_{L^2(\mathbf{R})}\int_0^t\|\bar u_x\|_{L^\i(\mathbf{R})}d\tau\\
&\di \leq
C\sqrt\v\ln(1+T)\Big[\sup_{t\in[0,T]}\|\sqrt\v(\r^{\t-\f12})_x\|^2_{L^2(\mathbf{R})}+\sup_{t\in[0,T]}\|\sqrt\r(u-\bar
u)\|^2_{L^2(\mathbf{R})}\Big].
\end{array}
$$
Thus the term $I_2$ is estimated as
\begin{eqnarray}\label{11-12-1}
\int_0^t\int_{\mathbf{R}}I_2dxd\tau &&\le
C\sqrt\v\ln(1+T)\Big[1+\sup_{t\in[0,T]}\|\sqrt\v(\r^{\t-\f12})_x\|^2_{L^2(\mathbf{R})}+\sup_{t\in[0,T]}\|\sqrt\r(u-\bar
u)\|^2_{L^2(\mathbf{R})}\Big]\nonumber\\
&&+C\v \sup_{t\in [0,T]}\|\r\Psi(\r,\bar\r)\|_{L^1(\mathbf{R})}.
\end{eqnarray}
The term $I_3$ can be estimated as follows.
\begin{equation}
\begin{array}{ll}
\di
\int_0^t\int_{\mathbf{R}}I_3dxd\tau&\di=\int_0^t\int_{\mathbf{R}}I_3({\bf
1}_{\{|\r-\bar\r|\leq \f{\r_{-}} 2\}}+{\bf 1}_{\{|\r-\bar\r|>\f{\r_{-}} 2\}})dxd\tau\\
&\di:=J_3^1+J_3^2.
\end{array}
\end{equation}
Direct estimates lead to
\begin{equation}
\begin{array}{ll}
J_3^1&\di=\int_0^t\int_{\mathbf{R}}\f{8\a\g}{(\a+\g-1)^2}(\bar\r^{\f{\a+\g-1}{2}})_{xx}(\r^{\f{\a+\g-1}{2}}-\bar\r^{\f{\a+\g-1}{2}}){\bf
1}_{\{|\r-\bar\r|\leq \f{\r_{-}} 2\}}dxd\tau\\[3mm]
&\di\le
C\sup_{t\in[0,T]}\|(\r^{\f{\a+\g-1}{2}}-\bar\r^{\f{\a+\g-1}{2}}){\bf
1}_{\{|\r-\bar\r|\leq \f{\r_{-}} 2\}}\|_{L^2(\mb{R})}\int_0^t\|(\bar\r^{\f{\a+\g-1}{2}})_{xx}\|_{L^2(\mathbf{R})}d\tau\\[3mm]
&\leq C\|\sqrt{\r\Psi(\r,\bar\r)}\|_{L^2(\mathbf{R})},
\end{array}
\end{equation}
where one has used the fact that
$$
\|(\r^{\f{\a+\g-1}{2}}-\bar\r^{\f{\a+\g-1}{2}}){\bf
1}_{\{|\r-\bar\r|\leq \f{\r_{-}} 2\}}\|_{L^2(\mb{R})}\le C
\|\sqrt{\r\Psi(\r,\bar\r)}\|_{L^2(\mathbf{R})}.
$$

 Moreover, due to the facts that
$$
\lim_{\r\rightarrow0+}\f{|\r^{\f{\a+\g-1}{2}}-\bar\r^{\f{\a+\g-1}{2}}|{\bf
1}_{\{|\r-\bar\r|>\f{\bar\r}
2\}}}{\r\Psi(\r,\bar\r)}=\bar\r^{\f{\a-\g-1}{2}}\leq C,
$$
and for
\begin{equation}
\f{\a+\g-1}{2}\leq \g, \quad{\rm i.\ e.},\quad \a\leq
\g+1,\label{(3.36)}
\end{equation}
$$
\lim_{\r\rightarrow+\i}\f{|\r^{\f{\a+\g-1}{2}}-\bar\r^{\f{\a+\g-1}{2}}|{\bf
1}_{\{|\r-\bar\r|>\f{\r_{-}} 2\}}}{\r\Psi(\r,\bar\r)}\leq C,
$$
we can estimate $J_3^2$ as
\begin{equation}
\begin{array}{ll}
\di J_3^2&\di
=\int_0^t\int_{\mathbf{R}}\f{8\a\g}{(\a+\g-1)^2}(\bar\r^{\f{\a+\g-1}{2}})_{xx}(\r^{\f{\a+\g-1}{2}}-\bar\r^{\f{\a+\g-1}{2}}){\bf
1}_{\{|\r-\bar\r|>\f{\r_{-}} 2\}}dxd\tau\\
&\di\le
C\sup_{t\in[0,T]}\|(\r^{\f{\a+\g-1}{2}}-\bar\r^{\f{\a+\g-1}{2}}){\bf 1}_{\{|\r-\bar\r|>\f{\r_{-}} 2\}}\|_{L^1(\mb{R})}\int_0^t\|(\bar\r^{\f{\a+\g-1}{2}})_{xx}\|_{L^\i(\mathbf{R})}d\tau\\[3mm]
&\leq C\eta^\f{2}{4q+1}\|\r\Psi(\r,\bar\r)\|_{L^1(\mathbf{R})}.
\end{array}
\end{equation}

The term $I_4$ can be handled similarly because $0<\t<\f12<\a\leq
\g+1$. Now we turn to
$$
\begin{array}{ll}
\di \int_0^t\int_{\mathbf{R}}I_5 dxd\tau &\di
=\int_0^t\int_{\mathbf{R}}I_5({\bf
1}_{\{|\r-\bar\r|\leq \f{\r_{-}} 2\}}+{\bf 1}_{\{|\r-\bar\r|>\f{\r_{-}} 2\}}) dxd\tau\\
&\di :=J_5^1+J_5^2.
\end{array}
$$
One has
$$
\begin{array}{ll}
\di J_5^1&\di
=\int_0^t\int_{\mathbf{R}}\f{2\g}{(\a+\g-1)}\Big[(\bar\r^{\f{\a+\g-1}{2}})_x\bar\r^{\f{\g-\a-1}{2}}\Big]_x(\r^\a-\bar\r^\a){\bf
1}_{\{|\r-\bar\r|\leq \f{\r_{-}} 2\}}dxd\tau\\
&\di \leq
C\sup_{t\in[0,T]}\|(\r^\a-\bar\r^\a){\bf1}_{\{|\r-\bar\r|\leq
\f{\r_{-}} 2\}}\|_{L^2(\mb{R})}\int_0^t\|\Big[(\bar\r^{\f{\a+\g-1}{2}})_x\bar\r^{\f{\g-\a-1}{2}}\Big]_x\|_{L^2(\mb{R})}d\tau\\
&\di \leq
C\sup_{t\in[0,T]}\|\sqrt{\r\Psi(\r,\bar\r)}\|_{L^2(\mathbf{R})},
\end{array}
$$
and
$$
\begin{array}{ll}
J_5^2
&\di=\int_0^t\int_{\mathbf{R}}\f{2\g}{(\a+\g-1)}\Big[(\bar\r^{\f{\a+\g-1}{2}})_x\bar\r^{\f{\g-\a-1}{2}}\Big]_x(\r^\a-\bar\r^\a)
{\bf 1}_{\{|\r-\bar\r|>\f{\r_{-}} 2\}}dxd\tau\\
&\di \leq C\sup_{t\in[0,T]}\|(\r^\a-\bar\r^\a){\bf1}_{\{|\r-\bar\r|>
\f{\r_{-}} 2\}}\|_{L^1(\mb{R})}\int_0^t\|\Big[(\bar\r^{\f{\a+\g-1}{2}})_x\bar\r^{\f{\g-\a-1}{2}}\Big]_x\|_{L^\i(\mb{R})}d\tau\\
&\di \leq
C\eta^{\f2{4q+1}}\sup_{t\in[0,T]}\|\sqrt{\r\Psi(\r,\bar\r)}\|_{L^1(\mathbf{R})},
\end{array}
$$
where the following facts have been used:
$$
\lim_{\r\rightarrow0+}\f{|\r^\a-\bar\r^\a|{\bf
1}_{\{|\r-\bar\r|>\f{\r_{-}}
2\}}}{\r\Psi(\r,\bar\r)}=\bar\r^{\a-\g}\leq C,
$$
\begin{equation}
\lim_{\r\rightarrow+\i}\f{|\r^\a-\bar\r^\a|{\bf
1}_{\{|\r-\bar\r|>\f{\r_{-}} 2\}}}{\r\Psi(\r,\bar\r)}\leq C, \quad
{\rm for} \quad \f12<\a\leq \g.\label{(3.38)}
\end{equation}
Finally, $I_6$ can be estimated as for $I_5$. In fact, we note that
the term $I_5$ involves the index $\alpha$ and $I_6$ involves the
index $\theta$. Since $0<\theta<\frac12<\alpha$, when we make the
estimate in the case $\rho\to+\infty$(see \eqref{(3.38)}), the order
of $\alpha$ will be dominant and hence the estimate of $I_6$ is much
more direct.

Now for $\a$ and $\g$ satisfying \eqref{(3.2)}, we can obtain
\begin{equation}
\begin{array}{ll}
\di \sup_{t\in[0,T]}\int_{\mathbf{R}}\Big\{\r\left[(u-\bar
u)+(\varphi_\v^{\a,\t}(\r))_x\right]^2+\r(u-\bar
u)^2+\r\Psi(\r,\bar\r)\Big\}(x,t)dx\\
\di +\int_0^T\int_{\mathbf{R}}\Big\{\bar
u_x\Big[p(\r)-p(\bar\r)-p^\prime(\bar\r)(\r-\bar\r)\Big]+\r(u-\bar
u)^2\bar u_x+(\r^\a+\v\r^\t)\Big[(u-\bar u)_x\Big]^2\\
\di +\Big[(\r^{\f{\a+\g-1}{2}}-\bar\r^{\f{\a+\g-1}{2}})_x\Big]^2
+\v\Big[(\r^{\f{\t+\g-1}{2}}-\bar\r^{\f{\t+\g-1}{2}})_x\Big]^2\Big\}(x,t)dxdt\\
\di \leq C+C\sqrt\v
\ln(1+T)\Big\{1+\sup_{t\in[0,T]}\int_{\mathbf{R}}\Big[\r(u-\bar
u)^2+\r(\varphi_\v^{\a,\t}(\r))_x^2\Big]dx\Big\}\\
\di \leq C+C\sqrt\v
\ln(1+T)\Big\{1+\sup_{t\in[0,T]}\int_{\mathbf{R}}\r(u-\bar
u)^2+\r\Big[(u-\bar u)+(\varphi_\v^{\a,\t}(\r))_x\Big]^2dx\Big\}
\end{array}
\end{equation}
where $C>0$ is the constant independent of $\v$ and $t$.

Choosing $\v$ such that $\sqrt\v \ln(1+T)\leq \v^{\f14}$ and $\v$
small enough, we arrive at
\begin{equation}
\begin{array}{ll}
\di \sup_{t\in[0,T]}\int_{\mathbf{R}}\Big\{\r\left[(u-\bar
u)+(\varphi_\v^{\a,\t}(\r))_x\right]^2+\r(u-\bar
u)^2+\r\Psi(\r,\bar\r)\Big\}(x,t)dx\\
\di +\int_0^T\int_{\mathbf{R}}\Big\{\bar
u_x\Big[p(\r)-p(\bar\r)-p^\prime(\bar\r)(\r-\bar\r)\Big]+\r(u-\bar
u)^2\bar u_x+(\r^\a+\v\r^\t)\Big[(u-\bar u)_x\Big]^2\\
\di +\Big[(\r^{\f{\a+\g-1}{2}}-\bar\r^{\f{\a+\g-1}{2}})_x\Big]^2
+\v\Big[(\r^{\f{\t+\g-1}{2}}-\bar\r^{\f{\t+\g-1}{2}})_x\Big]^2\Big\}(x,t)dxdt\\
\di \leq C
\end{array}
\end{equation}
Consequently, combining \eqref{(3.29)}, \eqref{(3.36)} and
\eqref{(3.38)} shows that for $\a$ and $\g$ satisfying
\eqref{(3.2)}, it holds that
\begin{equation}
\begin{array}{ll}
\di \sup_{t\in[0,T]}\int_{\mathbf{R}}\Big\{\r(u-\bar
u)^2+\Big[\Big(\f{\r^{\a-\f12}}{\a-\f12}\Big)_x\Big]^2+\v^2\Big[\Big(\f{\r^{\t-\f12}}{\t-\f12}\Big)_x\Big]^2+\r\Psi(\r,\bar\r)\Big\}(x,t)dx\\
\di +\int_0^T\int_{\mathbf{R}}\Big\{\bar
u_x\Big[p(\r)-p(\bar\r)-p^\prime(\bar\r)(\r-\bar\r)\Big]+\r(u-\bar
u)^2\bar u_x+(\r^\a+\v\r^\t)\Big[(u-\bar u)_x\Big]^2\\
\di +\Big[(\r^{\f{\a+\g-1}{2}}-\bar\r^{\f{\a+\g-1}{2}})_x\Big]^2
+\v\Big[(\r^{\f{\t+\g-1}{2}}-\bar\r^{\f{\t+\g-1}{2}})_x\Big]^2\Big\}(x,t)dxdt\\
\di \leq C.
\end{array}
\end{equation}
Thus lemma 3.1 is proved.

The following lemma is the key point to get the existence of the
approximate solution $(\r_\v,u_\v)(t,x)$.

{\em Lemma 3.2} There exist an absolutely constant $C$ and a
positive constant $C(\v,T)$ depending on $\v$ and $T$  such that
\begin{eqnarray}
0<C(\v,T)\le\rho_\v\le C.\label{D1}
\end{eqnarray}

{\em Proof:} First, we derive the upper bound for $\r_\v(x,t)$.

It follows from the entropy estimate that
\begin{eqnarray}
&&(\rho^{\alpha-\frac12}_\v(x,t)-\bar\rho^{\alpha-\frac12}(x,t))^{2}
=\int_{-\i}^x
[(\rho^{\alpha-\frac12}_\v-\bar\rho^{\alpha-\frac12})^{2}]_x
dx\nonumber\\
&&=2\int_{-\i}^x
(\rho^{\alpha-\frac12}_\v-\bar\rho^{\alpha-\frac12})(\rho_\v^{\alpha-\frac12}-\bar\rho^{\alpha-\frac12})_x
dx\nonumber\\
&&\le 2(\int_{-\i}^x
|(\rho^{\alpha-\frac12}_\v-\bar\rho^{\alpha-\frac12})|^{2}
dx)^\frac12(\int_{-\i}^x
[(\rho_\v^{\alpha-\frac12}-\bar\rho^{\alpha-\frac12})_x]^2
dx)^\frac12\nonumber\\
&&\le C+\int_{-\i}^x
|(\rho^{\alpha-\frac12}_\v-\bar\rho^{\alpha-\frac12})|^{2}[{\bf
1}|_{\{|\rho_\v-\bar\rho|<\frac{\bar\rho}{2}\}}+{\bf
1}|_{\{|\rho_\v-\bar\rho|\ge\frac{\bar\rho}{2}\}}] dx\nonumber\\
&&\equiv C+I_1(t)+I_2(t),\label{37-4}
\end{eqnarray}
for any fixed $t\in [0,T]$. Note that for
$|\rho_\v-\bar\rho|<\frac{\bar\rho}{2}$, that is,
$\frac{\bar\rho}{2}<|\rho_\v|<\frac32\bar\rho$, one has
\begin{eqnarray*}
|\rho_\v^{\alpha-\frac12}-\bar\rho^{\alpha-\frac12}|^2\le
C|\rho_\v-\bar\rho|^2\leq C\r_\v\Psi(\r_\v,\bar\r).
\end{eqnarray*}
Hence,
\begin{eqnarray}
I_1(t)\le C\|\r_\v\Psi(\r_\v,\bar\r)\|_{L^1(\mb{R})}\leq
C.\label{37-5}
\end{eqnarray}
If $|\rho_\v-\bar\rho|\ge\frac{\bar\rho}{2}$, then
$$
\lim_{\r_\v\rightarrow0+}\f{|\rho_\v^{\alpha-\frac12}-\bar\rho^{\alpha-\frac12}|^2{\bf
1}|_{\{|\rho_\v-\bar\rho|\ge\frac{\bar\rho}{2}\}}}{\r_\v\Psi(\r_\v,\bar\r)}=\bar\r^{2\a-1-\g}\leq
C,
$$
and
$$
\lim_{\r_\v\rightarrow+\i}\f{|\rho_\v^{\alpha-\frac12}-\bar\rho^{\alpha-\frac12}|^2{\bf
1}|_{\{|\rho_\v-\bar\rho|\ge\frac{\bar\rho}{2}\}}}{\r_\v\Psi(\r_\v,\bar\r)}\leq
C,
$$
if $\f12<\a\leq \f{\g+1}{2}$.

Hence
\begin{eqnarray}
I_2(t)\le C\|\r_\v\Psi(\r_\v,\bar\r)\|_{L^1(\mathbf{R})} \leq C.
\label{37-8}
\end{eqnarray}
It follows from \eqref{37-4}, \eqref{37-5}, \eqref{37-8} that
\begin{eqnarray}
(\rho_\v^{\alpha-\frac12}-\bar\rho^{\alpha-\frac12})^{2}\leq C,
\end{eqnarray}
which implies that
\begin{eqnarray}
|\rho_\v|\le C.
\end{eqnarray}
The upper bound of the approximate solution $\r_\v(t,x)$ is proved.

Next we derive a lower bound for $\r_\v(t,x)$. Since
$\lim_{\r\rightarrow0}\r\Psi(\r,\bar\r)=\bar\r^{\g}$, then
$\r_\v\Psi(\r_\v,\bar\r)$ is bounded away from 0 on
$[0,\f12\bar\r]$. Thus we can deduce from the bound on
$\r_\v\Psi(\r_\v,\bar\r)$ in $L^\i(0,T;L^1(\mathbf{R}))$ that
there exists a constant $C_1=C_1(T)>0$, such that for all
$t\in[0,T]$,
$$
{\rm meas} \{x\in\mathbf{R}|\r_\v(x,t)\leq \f12\bar\r(x,t)\}\leq
\f{1}{\inf_{\r\in[0,\f12\bar\r]}\r\Psi(\r,\bar\r)}\int_{\{x\in\mathbf{R}|\r_\v(x,t)\leq
\f12\bar\r(x,t)\}}\r_\v\Psi(\r_\v,\bar\r)(x,t)dx\leq C_1.
$$

Therefore, for every $x_0\in\mathbf{R}$, there exists $M=M(T)>0$
large enough, such that
$$
\begin{array}{ll}
\di \int_{|x-x_0|\leq M}\r_\v(x,t)dx&\di \geq\int_{\{|x-x_0|\leq
M\}\cap\{x\in\mathbf{R}|\r_\v(x,t)> \f12\bar\r(x,t)\}}\r_\v(x,t)dx\\
&\di \geq \f12\inf_{(x,t)}\bar\r(x,t){\rm meas}\Big\{\{|x-x_0|\leq
M\}\cap\{x\in\mathbf{R}|\r_\v(x,t)> \f12\bar\r(x,t)\}\Big\}\\
&\di =\f12\r_-{\rm meas}\Big\{\{|x-x_0|\leq
M\}\cap\{x\in\mathbf{R}|\r_\v(x,t)\leq \f12\bar\r(x,t)\}^{c}\Big\}\\
&\di \geq \f12 \r_-(2M-C_1)>0,
\end{array}
$$
for all $t\in [0,T]$.

From the continuity of $\r_\v$, there exists $x_1\in [x_0-M,x_0+M]$
such that
$$
\r_\v(x_1,t)=\int_{|x-x_0|\leq M}\r_\v(x,t)dx\geq \f12 \r_-(2M-C_1).
$$
Thus,
$$
\begin{array}{ll}
\di \r_\v^{\t-\f12}(x_0,t)&\di =\r_\v^{\t-\f12}(x_1,t)+\int_{x_1}^{x_0}(\r_\v^{\t-\f12})_x(x,t)dx\\
&\di \leq [\f12
\r_-(2M-C_1)]^{\t-\f12}+\|(\r_\v^{\t-\f12})_x(\cdot,t)\|_{L^2(\mathbf{R})}|x_1-x_0|^{\f12}\\
&\di \leq [\f12 \r_-(2M-C_1)]^{\t-\f12}+C_\v M^{\f12},
\end{array}
$$
where we have used the fact $0<\t<\f12$. Consequently, we can get
that
$$
\r_\v(x_0,t)\geq \Big\{[\f12 \r_-(2M-C_1)]^{\t-\f12}+C_\v
M^{\f12}\Big\}^{\f{2}{2\t-1}}:=C(\v,T).
$$
for any $x_0\in\mathbf{R}$ and $t\in[0,T]$.

With the lower and upper bounds on $\r_\v$, we can get the
existence of the approximate solution $(\r_\v,u_\v)(t,x)$ by a
similar argument as in \cite{MV2}. In order to pass the limit
$\v\ra0$, we need the following higher estimates on the momentum.

{\em Lemma 3.3} There exists a positive constant $C(T)$
independent of $\v$, such that
$$
\sup_{t\in[0,T]}\int_{\mathbf{R}}\r_\v|u_\v-\bar
u|^{3}(t,x)dx+\int_0^T\int_{\mathbf{R}}(\r_\v^\a+\v\r_\v^\t)[(u_\v-\bar
u)_x]^2|u_\v-\bar u| dxdt\leq C(T).
$$

{\em Proof:} Multiplying by $(u-\bar u)|u-\bar u|$ on the both
sides of \eqref{(3.5)} yields
\begin{equation}
\begin{array}{ll}
\di \big(\f23\r|u-\bar u|^3\big)_t+\big(\f23\r u|u-\bar u|^3\big)_x
+(p(\r)-p(\bar\r))_x(u-\bar u)|u-\bar u|\\\di +[(\r-\bar\r)\bar
u_t+(\r u-\bar\r\bar u)\bar u_x](u-\bar u)|u-\bar u|\\
\di =[(\r^\a+\v\r^\t)u_x(u-\bar u)|u-\bar
u|]_x-(\r^\a+\v\r^\t)[(u-\bar u)_x]^2|u-\bar u|\\
-(\r^\a+\v\r^\t)\bar u_x(u-\bar u)_x|u-\bar
u|-(\r^\a+\v\r^\t)u_x(u-\bar u)(|u-\bar u|)_x.
\end{array}
\end{equation}
Note that
$$
(|u-\bar u|)_x=sgn(u-\bar u)(u-\bar u)_x.
$$

Thus
\begin{equation}
\begin{array}{ll}
\di \big(\f23\r|u-\bar
u|^3\big)_t+H_{3x}(t,x)+2(\r^\a+\v\r^\t)[(u-\bar u)_x]^2|u-\bar
u|\\[3mm]
\di =-(u-\bar u)|u-\bar u|\Big[(\r-\bar\r)\bar u_t+(\r u-\bar\r\bar
u)\bar
u_x+(p(\r)-p(\bar\r))_x\Big]\\
\di \quad -2(\r^\a+\v\r^\t)\bar u_x(u-\bar u)_x|u-\bar
u|,\\
\di\quad
\end{array}\label{(3.55-)}
\end{equation}
where
$$
H_3(t,x)=\f23\r u|u-\bar u|^3-(\r^\a+\v\r^\t)u_x(u-\bar u)|u-\bar
u|.
$$
This, together with \eqref{(3.16)}, implies
\begin{equation}
\begin{array}{ll}
\di\big(\f23\r|u-\bar
u|^3\big)_t+H_{3x}(t,x)+2(\r^\a+\v\r^\t)[(u-\bar u)_x]^2|u-\bar
u|+\r|u-\bar u|^{3}\bar u_x \\[3mm]
\di =-p(\r)_x(u-\bar u)|u-\bar u|+\r\f{p(\bar\r)_x}{\bar\r}(u-\bar u)|u-\bar u|-2(\r^\a+\v\r^\t)\bar u_x(u-\bar u)_x|u-\bar u|\\
\di :=I_7+I_8+I_9.
\end{array}\label{(3.55)}
\end{equation}
Integrating \eqref{(3.55)} over $[0,t]\times\mathbf{R}$ with respect
to $t,x$ gives
\begin{equation}
\begin{array}{ll}
\di \int_{\mathbf{R}}\f23\r|u-\bar
u|^{3}(t,x)dx+\int_0^t\int_{\mathbf{R}}\Big[2(\r^\a+\v\r^\t)[(u-\bar
u)_x]^2|u-\bar
u|^\d+\r|u-\bar u|^{3}\bar u_x\Big] dxd\tau\\[3mm]
\di =\int_0^t\int_{\mathbf{R}} (I_7+I_8+I_9)dxd\tau.
\end{array} \label{(3.57)}
\end{equation}
Now we estimate such term on the right hand side of
\eqref{(3.57)}. First, integrating by parts with respect to $x$
gives
\begin{equation}
\begin{array}{ll}
\di \int_0^t\int_{\mathbf{R}} I_7dxd\tau=2\int_0^t\int_{\mathbf{R}}
p(\r)(u-\bar u)_x|u-\bar u|dxd\tau\\
\di \leq \int_0^t\int_{\mathbf{R}}\r^\a[(u-\bar
u)_x]^2dxd\tau+\int_0^t\int_{\mathbf{R}}\r^{2\g-\a-1}\r|u-\bar u|^2
dxd\tau\\
\di \leq C
\end{array}
\label{(3.58)}
\end{equation}
uniformly with respect to $\v$ if
\begin{equation}
2\g-\a-1\geq 0,\quad {\rm i.\ e.}\quad \a\leq 2\g-1.\label{(3.53)}
\end{equation}
In fact, the relation \eqref{(3.3)} guarantees \eqref{(3.53)}
because $\f{\g+1}{2}<2\g-1.$

Next,
\begin{equation}
\begin{array}{ll}
\di  \int_0^t\int_{\mathbf{R}} I_8dxd\tau&\di =
\int_0^t\int_{\mathbf{R}}\r\f{p(\bar\r)_x}{\bar\r}(u-\bar u)|u-\bar u| dxd\tau\\
&\di \leq C\sup_{t\in[0,T]}\int_{\mathbf{R}}\r|u-\bar
u|^2(t,x)dx\int_0^T\|\bar\r_x\|_{L^\i(\mathbf{R})}dt
\\
&\di\leq  C\sup_{t\in[0,T]}\int_{\mathbf{R}}\r|u-\bar
u|^2(t,x)dx\int_0^T(1+t)^{-1}dt \\
&\di \leq C_1\ln(1+T).
 \end{array}
\label{(3.59)}
\end{equation}

Finally,
\begin{equation}
\begin{array}{ll}
 \di\int_0^t\int_{\mathbf{R}} I_9dxd\tau&\di =
\int_0^t\int_{\mathbf{R}}-2(\r^\a+\v\r^\t)\bar u_x(u-\bar
u)_x|u-\bar u|dxd\tau\\
&\di \leq \int_0^t \int_{\mathbf{R}}(\r^\a+\v\r^\t)[(u-\bar
u)_x]^2|u-\bar u|dxd\tau+\int_0^t
\int_{\mathbf{R}}(\r^\a+\v\r^\t)\bar u_x^2|u-\bar u|dxd\tau.
\end{array}
\label{(3.60)}
\end{equation}
One also has
\begin{equation}
\begin{array}{ll}
\di \int_0^t \int_{\mathbf{R}}(\r^\a+\v\r^\t)\bar u_x^2|u-\bar
u|dxd\tau\\
\di\leq\int_0^t \int_{\mathbf{R}}\r|u-\bar u|^3dxd\tau+C\int_0^t
\int_{\mathbf{R}}(\r^{\f{3\a-1}{2}}+\v^{\f32}\r^{\f{3\t-1}{2}})\bar
u_x^3dxd\tau\\
\di \leq\int_0^t \int_{\mathbf{R}}\r|u-\bar u|^3dxd\tau+C\int_0^t
\|\bar u_x\|^3_{L^3(\mathbf{R})}d\tau
\\
\di \leq\int_0^t \int_{\mathbf{R}}\r|u-\bar
u|^3dxd\tau+C\int_0^t(1+\tau)^{-2}d\tau\\
\di \leq\int_0^t \int_{\mathbf{R}}\r|u-\bar u|^3dxd\tau+C
\end{array}
\label{(3.61)}
\end{equation}
if $\a,\t\geq\f13$. Without loss of generality, we can set
$\t=\f13.$

Substituting \eqref{(3.58)}-\eqref{(3.61)} into \eqref{(3.57)}
implies Lemma 3.3.

Now with these uniform in $\v$ estimates in hand, we can pass the
limit process $\v\rightarrow0$, obtain the existence of the weak
solution $(\r,u)(t,x)$, and get the uniform in time estimates in
Theorem 2.1.

\section{Asymptotic behavior of weak solutions}

In this section, we will study the asymptotic behavior of the weak
solution $(\r,u)(t,x)$ obtained in previous section. We assume
that the solution is smooth enough. The rigorous proof can be
obtained by using the usual regularization procedure.

{\bf Proof of Theorem 2.2.} Since $0\leq \rho\le C,
0<\r_{-}<\bar\r<\r_{+}$, it holds that
\begin{equation}
C_1^{-1}(\r-\bar\r)^2\leq \r\Psi(\r,\bar\r)\leq C_1(\r-\bar\r)^2
\end{equation}
for some constant $C_1>0$ which may depend on $C, \r_{-}, \r_{+}$.

In the following, we denote by $C>0$ a universal constant. For any
$s\ge 1$, Lemma 3.2 implies
\begin{eqnarray*}
|\rho^s-\bar\rho^s|^2\le C|\rho-\bar\rho|^2.
\end{eqnarray*}
Hence
\begin{eqnarray}
\int_{\mathbf{R}} |\rho^s-\bar\rho^s|^2 dx \le C\int_{\mathbf{R}}
|\rho-\bar\rho|^2 dx\le C.\label{Jan6-2}
\end{eqnarray}
Similarly,
\begin{eqnarray}
\int_{\mathbf{R}} |\rho^s-\bar\rho^s|^{2\l} dx \le
C\int_{\mathbf{R}} |\rho-\bar\rho|^{2\l} dx\le C\label{Jan6-4}
\end{eqnarray}
for any $\lambda\ge 1$. Moreover, one has
\begin{eqnarray*}
&&\int_{\mathbf{R}} |[(\rho^s-\bar\rho^s)^{2\l}]_x| dx=2\lambda
s\int_{\mathbf{R}}
|(\rho^s-\bar\rho^s)^{2\lambda-1}[\rho^{s-1}\rho_x-\bar\r^{s-1}\bar\r_x]|
dx\\
&&\le \frac{2\lambda s(2\alpha-1)}{2}(\int_{\mathbf{R}}
(\rho^s-\bar\rho^s)^{2(2\lambda-1)}\rho^{2s+1-2\alpha} dx)^\frac12
(\int_{\mathbf{R}} [(\rho^{\alpha-\frac12})_x]^2 dx)^\frac12\\
&&\quad +2\l s(\int_{\mathbf{R}}
|\rho^s-\bar\rho^s|^{2(2\lambda-1)}dx)^{\f12}(\int_{\mathbf{R}}
\bar\r^{2(s-1)}|\bar\r_x|^2dx)^{\f12}\\&&\leq C.
\end{eqnarray*}
It follows from \eqref{Jan6-4} that for any fixed $t$,
\begin{eqnarray}
\rho^s-\bar\rho^s\rightarrow 0
\end{eqnarray}
as $|x|\to \infty$. By  Lemma 3.1, it holds that
\begin{eqnarray}
\int_0^t\int_{\mathbf{R}}[(\rho^{\frac{\gamma+\alpha-1}{2}}-\bar\rho^{\frac{\gamma+\alpha-1}{2}})_x]^2
dxdt\le C,
\end{eqnarray}
with $C$ an absolute constant depending only on the initial data.
Set $b=\frac{\gamma+\alpha-1}{2}$. Then
\begin{eqnarray}
\int_0^t\int_{\mathbf{R}} [(\rho^b-\bar\r^b)_x]^2 dxdt\le
C.\label{Jan6-3}
\end{eqnarray}
Choosing $s>b+1$, one has
\begin{eqnarray*}
&&(\rho^s-\bar\rho^s)^2(t,x)=\int_{-\infty}^x
[(\rho^s-\bar\rho^s)^2]_x
dx=2\int_{-\infty}^x (\rho^s-\bar\rho^s)(\rho^s-\bar\r^s)_x dx\\
&&=2s\int_{-\infty}^x (\rho^s-\bar\rho^s)
(\rho^{s-1}\rho_x-\bar\r^{s-1}\bar\r_x)
dx\\
&&=\frac{2s}{b}\int_{-\infty}^x
(\rho^s-\bar\rho^s)[(\rho^b-\bar\r^b)_x\rho^{s-b}+(\bar\r^b)_x(\r^{s-b}-\bar\r^{s-b})] dx\\
&&\le C
\|\rho^s-\bar\rho^s\|_{L^2(\mathbf{R})}\|(\rho^b-\bar\r^b)_x\|_{L^2(\mathbf{R})}+C\int_\mathbf{R}(\bar\r^b)_x(\r^s-\bar\r^s)(\r^{s-b}-\bar\r^{s-b})dx\\
&&\leq
\|\rho^s-\bar\rho^s\|_{L^2(\mathbf{R})}\|(\rho^b-\bar\r^b)_x\|_{L^2(\mathbf{R})}+C\int_\mathbf{R}\bar
u_x[\r^\g-\bar\r^\g-\g\bar\r^{\g-1}(\r-\bar\r)]dx
\end{eqnarray*}
where in the last inequality, we have used the fact that
\begin{equation}
\begin{array}{ll}
&\di (\bar\r^b)_x(\r^s-\bar\r^s)(\r^{s-b}-\bar\r^{s-b})\\
&\di =b\bar\r^{b-1}\bar\r_x(\r^s-\bar\r^s)(\r^{s-b}-\bar\r^{s-b})\\
&\di =\f{b\bar\r^{b}}{\sqrt{p^\prime(\bar\r)}}\bar
u_x(\r^s-\bar\r^s)(\r^{s-b}-\bar\r^{s-b})\\
&\di \leq C\bar u_x[\r^\g-\bar\r^\g-\g\bar\r^{\g-1}(\r-\bar\r)].
\end{array}
\end{equation}
Consequently,
\begin{eqnarray}
&&\int_0^t\sup_{x\in \mathbf{R}}(\rho^s-\bar\rho^s)^4 dt\le
C\sup_{t\in[0,T]}\|\rho^s-\bar\rho^s\|^2_{L^2(\mathbf{R})}\int_0^t
\|(\rho^b-\bar\r^b)_x\|_{L^2(\mathbf{R})}^2 dt\\
&&\qquad+C\sup_{t\in[,T]}\int_\mathbf{R}[\r^\g-\bar\r^\g-\g\bar\r^{\g-1}(\r-\bar\r)]dx\int_0^t\int_\mathbf{R}\bar
u_x[\r^\g-\bar\r^\g-\g\bar\r^{\g-1}(\r-\bar\r)]dxd\tau\\
&&\le C.
\end{eqnarray}
Moreover, applying \eqref{Jan6-4} leads to
\begin{eqnarray}
&&\int_0^t\int_\mathbf{R}
(\rho^s-\bar\rho^s)^4(\rho^s-\bar\rho^s)^{2l}
dxdt\nonumber\\
&&\le \int_0^t [\sup_{x\in
\mathbf{R}}(\rho^s-\bar\rho^s)^4\int_\mathbf{R}
(\rho^s-\bar\rho^s)^{2l} dx] dt\nonumber\\
&&\le \sup_t \int_\mathbf{R} (\rho^s-\bar\rho^s)^{2l} dx\int_0^t
\sup_{x\in \mathbf{R}}(\rho^s-\bar\rho^s)^4 dt\le C,
\end{eqnarray}
where $l\ge 1$ is any real number. Hence
\begin{eqnarray}\label{(4.12)}
\int_0^t\int_\mathbf{R} (\rho^s-\bar\rho^s)^{4+2l} dxdt\le C.
\end{eqnarray}
 Denote $f(t)=\int_{\mathbf{R}} (\rho^s-\bar\rho^s)^{4+2l} dx$. Then $f(t)\in
 L^1(0,\infty)\cap L^\infty(0,\infty)$ due to \eqref{Jan6-4} and \eqref{(4.12)}. Furthermore, direct
 calculations show that
\begin{equation}
\begin{array}{ll}
\di \frac{d}{dt} f(t)&\di =(4+2l)s\int_\mathbf{R}
(\rho^s-\bar\rho^s)^{3+2l}
(\rho^{s-1}\rho_t-\bar\r^{s-1}\bar\r_x) dx\\
&\di =-(4+2l)s\int_\mathbf{R}
(\rho^s-\bar\rho^s)^{3+2l}[\rho^{s-1}(\rho
u)_x-\bar\r^{s-1}(\bar\r\bar u)_x]
dx\\
&\di =(4+2l)(3+2l)s\int_\mathbf{R}
(\rho^s-\bar\rho^s)^{2+2l}(\rho^s-\bar\r^s)_x(\rho^{s-1}\rho
u-\bar\r^{s-1}\bar\r\bar u)
dx\\
&\di \quad+(4+2l)s(s-1)\int_\mathbf{R}
(\rho^s-\bar\rho^s)^{3+2l}(\rho^{s-2}\rho_x\rho
u-\bar\rho^{s-2}\bar\rho_x\bar\r \bar u)
dx\\
&\di =(4+2l)(3+2l)s\int_\mathbf{R}
(\rho^s-\bar\rho^s)^{2+2l}(\rho^s-\bar\r^s)_x(\rho^{s}u-\bar\r^{s}\bar
u)
dx\\
&\di \quad+(4+2l)s(s-1)\int_\mathbf{R}
(\rho^s-\bar\rho^s)^{3+2l}(\rho^{s-1}\rho_x
u-\bar\rho^{s-1}\bar\rho_x \bar u) dx\\
&\di =(4+2l)(3+2l)s\int_\mathbf{R}
(\rho^s-\bar\rho^s)^{2+2l}(\rho^s-\bar\r^s)_x\rho^{s}(u-\bar u)
dx\\
&\quad\di +(4+2l)(3+2l)s\int_\mathbf{R}
(\rho^s-\bar\rho^s)^{2+2l}(\rho^s-\bar\r^s)_x\bar
u(\rho^{s}-\bar\r^{s})
dx\\
&\di \quad+(4+2l)s(s-1)\int_\mathbf{R}
(\rho^s-\bar\rho^s)^{3+2l}\rho^{s-1}\rho_x
(u-\bar u) dx\\
&\di \quad+(4+2l)s(s-1)\int_\mathbf{R}
(\rho^s-\bar\rho^s)^{3+2l}\bar u(\rho^{s-1}\rho_x
-\bar\rho^{s-1}\bar\rho_x) dx\\
&\di :=J_1(t)+J_2(t)+J_3(t)+J_4(t).
\end{array}
\end{equation}
Now we claim that $J_i(t)\in L^2(0,+\i),~(i=1,2,3,4)$. In fact,
\begin{equation}
\begin{array}{ll}
\di J_1(t)&\di=\f{(4+2l)(3+2l)s^2}{b}\int_\mathbf{R}
(\rho^s-\bar\rho^s)^{2+2l}\rho^{s}(u-\bar
u)\Big[\r^{s-b}(\rho^b-\bar\r^b)_x+(\r^{s-b}-\bar\r^{s-b})(\bar\r^b)_x\Big]
dx\\
&\di =\f{(4+2l)(3+2l)s^2}{b}\int_\mathbf{R}
(\rho^s-\bar\rho^s)^{2+2l}\sqrt{\r}(u-\bar
u)\r^{2s-b-\f12}(\rho^b-\bar\r^b)_xdx\\
&\di\quad +\f{(4+2l)(3+2l)s^2}{b}\int_\mathbf{R}
(\rho^s-\bar\rho^s)^{2+2l}\sqrt{\r}(u-\bar
u)\r^{s-\f12}(\r^{s-b}-\bar\r^{s-b})(\bar\r^b)_x dx\\
&\leq C\|\sqrt{\r}(u-\bar
u)\|_{L^2(\mathbf{R})}\|(\rho^b-\bar\r^b)_x\|_{L^2(\mathbf{R})}\\
&\di\quad +C\|\sqrt{\r}(u-\bar
u)\|_{L^2(\mathbf{R})}\|(\rho^s-\bar\rho^s)^{2+2l}(\r^{s-b}-\bar\r^{s-b})(\bar\r^b)_x
\|_{L^2(\mathbf{R})}\\
&\leq C\|\sqrt{\r}(u-\bar
u)\|_{L^2(\mathbf{R})}\|(\rho^b-\bar\r^b)_x\|_{L^2(\mathbf{R})}\\
&\di\quad +C\|\sqrt{\r}(u-\bar u)\|_{L^2(\mathbf{R})}\|\bar
u_x[\r^\g-\bar\r^\g-\g\bar\r^{\g-1}(\r-\bar\r)] \|_{L^2(\mathbf{R})}
\end{array}
\end{equation}
Thus,
\begin{equation}
\begin{array}{ll}
\di \int_0^t |J_1(t)|^2dt&\di \leq
C\sup_{t\in[0,T]}\|\sqrt{\r}(u-\bar
u)\|^2_{L^2(\mathbf{R})}\int_0^t\|(\rho^b-\bar\r^b)_x\|^2_{L^2(\mathbf{R})}dt\\
&\di\quad +C\sup_{t\in[0,T]}\|\sqrt{\r}(u-\bar
u)\|^2_{L^2(\mathbf{R})}\int_0^t\|\bar
u_x[\r^\r-\bar\r^\g-\g\bar\r^{\g-1}(\r-\bar\r)]
\|_{L^2(\mathbf{R})}^2dt\\
&\di \leq C\sup_{t\in[0,T]}\|\sqrt{\r}(u-\bar
u)\|^2_{L^2(\mathbf{R})}\int_0^t\|(\rho^b-\bar\r^b)_x\|^2_{L^2(\mathbf{R})}dt\\
&\di\quad +C\sup_{t\in[0,T]}\|\sqrt{\r}(u-\bar
u)\|^2_{L^2(\mathbf{R})}\sup_{t\in[0,T]}\|\r\Psi(\r,\bar\r)\|_{L^1(\mathbf{R})}\cdot \\
&\di \quad ~\int_0^t\int_{\mathbf{R}}\bar
u_x[\r^\g-\bar\r^\g-\g\bar\r^{\g-1}(\r-\bar\r)] dxdt\\
&\di \leq C.
\end{array}
\end{equation}
The fact that $J_i(t)\in L^2(0,+\i),~(i=2,3,4)$ can be shown
similarly.

Hence
\begin{eqnarray}
\frac{d}{dt}f(t)\in L^2(0,+\infty).
\end{eqnarray}
Combining the obtained fact that $f(t)\in
 L^1(0,+\infty)\cap L^\infty(0,+\infty)$, one has
\begin{eqnarray}
f(t)\to 0, \ \ t\to+\infty.
\end{eqnarray}
Letting $m\ge 1$ be any real number to be determined later, we have
\begin{equation}
\begin{array}{ll}
\di |(\rho^s-\bar\rho^s)^m|&\di =|\int_{-\infty}^x
[(\rho^s-\bar\rho^s)^m]_x dx|\\
&\di =|m\int_{-\infty}^x
(\rho^s-\bar\rho^s)^{m-1}(\rho^s-\bar\r^s)_x dx|\\
&\di =|m\int_{-\infty}^x
(\rho^s-\bar\rho^s)^{m-1}\Big[\f{s}{\a-\f12}(\rho^{\a-\f12})_x\r^{s-\a+\f12}-s\bar\r^{s-1}\bar\r_x \Big]dx|\\
&\di \leq
C\|(\rho^s-\bar\rho^s)^{m-1}\|_{L^2(\mathbf{R})}\Big[\|(\rho^{\alpha-\frac12})_x\|_{L^2(\mathbf{R})}+\|\bar\r_x\|_{L^2(\mathbf{R})}\Big]\\
&\di \le C\|(\rho^s-\bar\rho^s)^{m-1}\|_{L^2(\mathbf{R})}.
\end{array}
\end{equation}
Choose $2(m-1)=4+2l$ to get
\begin{eqnarray}
\sup_{x\in \mathbf{R}}|(\rho^s-\bar\rho^s)^m|\le C f^\frac12(t)\to 0
\end{eqnarray}
as $t\to+\infty$.

Therefore,
$$
\lim_{t\to+\infty}\sup_{x\in\mathbf{R}}|\rho^s-\bar\rho^s|=0.
$$

 Now
we prove that
$\di\lim_{t\to+\infty}\sup_{x\in\mathbf{R}}|\rho-\bar\rho|=0$. Using
the fact that
$$
\begin{array}{ll}
\di |\rho-\bar\rho|^s&\di =|\rho-\bar\rho|^s{\bf
1}_{\{0\leq\rho\leq\f{\r_{-}}{2}\}}+|\rho-\bar\rho|^s{\bf
1}_{\{ \r>\f{\r_{-}}{2}\}}\\
&\di\leq C|\rho^s-\bar\rho^s|{\bf
1}_{\{0\leq\rho\leq\f{\r_{-}}{2}\}}+C|\rho^s-\bar\rho^s|^s{\bf
1}_{\{\rho>\f{\r_{-}}{2}\}}.
\end{array}
$$
Therefore, we have
\begin{eqnarray*}
&& \sup_{x\in \mathbf{R}}|\rho-\bar\rho|^s \leq
C\sup_{x\in\mathbf{R}}|\rho^s-\bar\rho^s|+C\sup_{x\in
\mathbf{R}}|\rho^s-\bar\rho^s|^s\to 0,
\end{eqnarray*}
 as $t\to+\infty$, which implies that
$$
 \lim_{t\to\infty}\sup_{x\in\mathbf{R}}|\rho-\bar\rho|=0.
$$
The proof of the lemma is finished.

\section{Vanishing of vacuum states and blow-up phenomena}

In this subsection, we first give a sketch of proof of Theorem 2.3
and then give some remarks on the blow-up phenomena of the solutions
when the vacuum states vanish. These results are similar to those in
\cite{LLX} in which the initial-boundary value problem and periodic
boundary value problem are studied.

{\em Proof:} [Sketch of proof of Theorem 2.3]

It follows from Theorem 2.2 that for any
$0<\rho_1<\inf_{t,x}\bar\rho(t,x)$, there exists a time $T_0>0$
such that
\begin{eqnarray}
0<\rho_1\le \rho(x,t)\le C, \ \ \ (x,t)\in \mathbf{R}\times
[T_0,\infty). \label{Jan22-1}
\end{eqnarray}
Therefore, for $t\ge T_0$, the density is bounded away from the
zero and the vacuum states vanish. Using the estimate in
\eqref{th12-3} and the standard theory for linear parabolic
equations, one can obtain that for $t\ge T_0$, the weak solution
becomes a strong solution to \eqref{(1.1)}-\eqref{(1.2)},
satisfying
\begin{eqnarray}\label{Jan22-2}
\left\{
\begin{array}{ll}
&\rho-\bar\rho\in L^\infty(T_0, t;H^1(\mathbf{R})), \ \ \rho_t\in
L^\infty(T_0,t; L^2(\mathbf{R})),\\
&u-\bar u\in L^2(T_0,t; H^2(\mathbf{R})), \ \ u_t\in L^2(T_0,t;
L^2(\mathbf{R})).
\end{array}
\right.
\end{eqnarray}
The detail of the proof is referred to \cite{LLX} and is omitted
it here. Furthermore, the asymptotic behaviors
$\lim_{t\to\infty}\sup_{x\in \mathbf{R}} |\rho-\bar\rho|=0$ and $
\lim_{t\to\infty} \|\rho-\bar\rho\|_{L^p}=0$ for $2<p\le \infty$
follow directly from \eqref{th14-1} and the estimate
$\|\rho-\bar\rho\|_{L^2}\le C$. The asymptotic behavior on the
velocity $\lim_{t\to\infty} \|u-\bar u\|_{L^2}=0$ follows from the
standard arguments, see \cite{SZ} for instance.

It should be remarked that we also have  finite blow-up phenomena
for the weak solutions of the Cauchy problem
\eqref{(1.1)}-\eqref{(1.3)} at the time when the vacuum states
vanish if the density contains vacuum states at least at one point.
These are similar as in \cite{LLX} in which the 1D initial-boundary
value problem and periodic problem are investigated. To be more
precise, we note that, if the density contains vacuum states at
least at one point, due to the facts that $\rho\in
C(\mathbf{R}\times [0,T])$ for any $T>0$ and
$\lim_{t\to\infty}\sup_{x\in \mathbf{R}} |\rho-\bar\rho|=0$, there
exists some critical time $T_1\in [0,T_0)$ with $T_0>0$ given by
\eqref{Jan22-1} and a nonempty subset $\Omega^0\subset \mathbf{R}$
such that
\begin{eqnarray}\label{Jan22-3}
\left\{
\begin{array}{ll}
&\rho(x,T_1)=0,  \ \forall x\in \Omega^0,\\
&\rho(x,T_1)>0,  \ \forall x\in \mathbf{R}\backslash\Omega^0,\\
&\rho(x,t)>0, \ \forall (x,t)\in \mathbf{R}\times (T_1,T_0].
\end{array}
\right.
\end{eqnarray}
From \eqref{Jan22-2}, it is easy to know that for any $\delta>0$,
$$
\int_{T_1+\delta}^{T_0} \|u_x\|_{L^\infty}
ds\leq\int_{T_1+\delta}^{T_0} \|(u-\bar u)_x\|_{L^\infty}
ds+\int_{T_1+\delta}^{T_0} \|\bar u_x\|_{L^\infty} ds <\infty.
$$
However, one has the following blow-up result of the solution.

{\em Theorem 5.1}\label{th5.1}
 Let $(\rho,u)$ be any global weak solution, which
contains vacuum states at least at one point for some time, to the
Cauchy problem \eqref{(1.1)}-\eqref{(1.3)} satisfying
\eqref{th12-2}-\eqref{th12-3}. Let $T_0>0$ and $T_1\in [0,T_0)$ be
the time such that \eqref{Jan22-1} and \eqref{Jan22-3} holds
respectively. Then, the solution $(\rho,u)$ blows up as vacuum
states vanish in the following sense: for any $\eta>0$, it holds
\begin{eqnarray}\label{Jan22-5}
\lim_{t\to T_1^+}\int_t^{T_1+\eta} \|u_x\|_{L^\infty} ds=\infty.
\end{eqnarray}

On the other hand, if there exists some $T_2\in (0,T_0)$ such that
the weak solution $(\rho,u)$ satisfies
$$
\|u-\bar u\|_{L^1(0,T_2;W^{1,\infty}(\mathbf{R}))}<\infty,
$$
then, there is a time $T_3\in [T_2, T_0)$ such that
\begin{eqnarray}\label{Jan22-6}
\lim_{t\to T_3^-}\int_0^t \|u_x\|_{L^\infty} ds=\infty.
\end{eqnarray}

The proof of Theorem 5.1 is completely similar to that in
\cite{LLX}. For completeness, we just sketch it here.

{\em Proof}
 It suffices to prove \eqref{Jan22-5} since the proof of
 \eqref{Jan22-6} is similar. If \eqref{Jan22-5} is not true, then
 there exists a fixed constant $\eta>0$, such that
\begin{eqnarray}\label{Jan22-7}
\int_{T_1}^{T_1+\eta} \|u_x\|_{L^\infty} ds<\infty.
\end{eqnarray}
Thanks to \eqref{Jan22-2} and \eqref{Jan22-7}, the particle path
$x(s)=X(s;t,x)$ through $(x,t)\in R\times (T_1,T_1+\eta]$ can be
well defined by solving
\begin{eqnarray}\label{Jan22-8}
\left\{
\begin{array}{lll}
&\frac{\partial}{\partial s}X(s;t,x)=u(X(s;t,x),s), \ \ &T_1\le
s<T_1+\eta,\\
 &X(t;t,x)=x, &T_1\le t<T_1+\eta, x\in R.
\end{array}
\right.
\end{eqnarray}
Then by the continuity equation \eqref{(1.1)}, one has
\begin{eqnarray}\label{Jan22-9}
\rho(x,t)=\rho(X(T_1;t,x),T_1)\exp\{-\int_{T_1}^t
u_y(y,s)|_{y=X(s;t,x)} ds\}
\end{eqnarray}
for any $(x,t)\in R\times (T_1,T_1+\eta]$. It follows from
\eqref{Jan22-7} and \eqref{Jan22-8} that  for $x_1\in \Omega_0$
defined by \eqref{Jan22-3}, which satisfies $\rho(x_1,T_1)=0$, there
exists a trajectory $x=x_1(t)\in R$ for $t\in [T_1,T_1+\eta]$ so
that $X(T_1;t,x_1(t))=x_1$. Thus, due to \eqref{Jan22-9} and
\eqref{Jan22-7}, one has that $\rho(x_1(t),t)=0$ for all $t\in
(T_1,T_1+\eta]$, which is a contradiction to \eqref{Jan22-3}.
\eqref{Jan22-5} is then proved and the proof of the theorem is
finished.

\end{document}